\documentclass[a4paper,11pt]{article}
\usepackage{amsmath,amsthm,amsfonts,amssymb,enumitem,graphicx,calc,microtype,mathtools,thmtools,underscore,anyfontsize,lmodern,float}
\usepackage{latexsym}
\usepackage[dvipsnames]{xcolor}
\usepackage{todonotes}
\usepackage{wrapfig}
\usepackage{paralist}
\usepackage{enumitem}
\usepackage[english]{babel}
\usepackage[utf8]{inputenc}
\usepackage[T1]{fontenc}
\usepackage[margin=30mm]{geometry}
\renewcommand{\baselinestretch}{1.1}
\allowdisplaybreaks
\newcommand{\defn}[1]{\textcolor{Maroon}{\emph{#1}}}
\usepackage[numbers,sort&compress,longnamesfirst]{natbib}
\makeatletter
\def\NAT@spacechar{~}
\makeatother
\usepackage[pdftitle={Shallow Minors, Graph Products and Beyond Planar Graphs},
pdfauthor={Hickingbotham -- Wood},
colorlinks =true,
linkcolor={black},
urlcolor={blue!80!black},
filecolor={blue!80!black},
citecolor={black},
anchorcolor=blue,
filecolor=blue, 
menucolor=blue]{hyperref}

\usepackage[noabbrev,capitalise]{cleveref}
\crefname{lem}{Lemma}{Lemmas}
\crefname{thm}{Theorem}{Theorems}
\crefname{cor}{Corollary}{Corollaries}
\crefname{prop}{Proposition}{Propositions}
\crefname{conj}{Conjecture}{Conjectures}
\crefname{open}{Open Problem}{Open Problems}
\crefname{obs}{Observation}{Observations}
\crefformat{equation}{(#2#1#3)}
\Crefformat{equation}{Equation #2(#1)#3}

\theoremstyle{plain}
\newtheorem{thm}{Theorem}
\newtheorem{lem}[thm]{Lemma}
\newtheorem{cor}[thm]{Corollary}
\newtheorem{obs}[thm]{Observation}
\newtheorem{prop}[thm]{Proposition}

\theoremstyle{definition}

\DeclarePairedDelimiter{\ceil}{\lceil}{\rceil}
\DeclarePairedDelimiter{\floor}{\lfloor}{\rfloor}

\DeclareMathOperator{\sn}{sn}
\DeclareMathOperator{\qn}{qn}
\DeclareMathOperator{\sqn}{sqn}

\DeclareMathOperator{\tw}{tw}

\renewcommand{\leq}{\leqslant}

\renewcommand{\geq}{\geqslant}
\newcommand{\subsetsim}{\mathrel{\substack{\textstyle\subset\\[-0.2ex]\textstyle\sim}}}
\newcommand{\beginproof}{\emph{Proof.}}

\theoremstyle{definition}

\DeclareMathOperator{\dist}{dist}

\DeclareMathOperator*{\ltw}{ltw}
\DeclareMathOperator*{\rtw}{rtw}

\DeclareMathOperator*{\bs}{\backslash}
\DeclareMathOperator*{\sse}{\subseteq}
\DeclareMathOperator{\R}{\mathbb{R}}
\DeclareMathOperator*{\N}{\mathbb{N}}

\DeclareMathOperator*{\G}{\mathcal{G}}
\DeclareMathOperator*{\Path}{\mathcal{P}}
\DeclareMathOperator*{\Pcal}{\mathcal{P}}

\DeclareMathOperator*{\B}{\mathcal{B}} 
 
\DeclareMathOperator*{\Tcal}{\mathcal{T}}

\DeclareMathOperator{\scol}{scol}
\DeclareMathOperator{\wcol}{wcol}
\DeclareMathOperator*{\cross}{cr}

\DeclareMathOperator*{\boxicity}{box}

\newcommand{\JournalArxiv}[2]{#1}

\presetkeys%
{todonotes}%
{inline}{}

\begin{document}
	
	\title{\bf\Large Shallow Minors, Graph Products\\ and Beyond Planar Graphs}
	\author{%
		Robert Hickingbotham\thanks{School of Mathematics, Monash University, Melbourne, Australia (\texttt{robert.hickingbotham@monash.edu}). Research supported by an Australian Government Research Training Program Scholarship.}
		\quad
		David R. Wood\thanks{School of Mathematics, Monash University, Melbourne, Australia (\texttt{david.wood@monash.edu}). Research supported by the Australian Research Council.}
	}
	
	\date{\normalsize\today}
	\maketitle
	
	\begin{abstract}
		 \JournalArxiv{The planar graph product structure theorem of Dujmovi\'{c}, Joret, Micek, Morin, Ueckerdt, and Wood [J. ACM 2020] states that every planar graph is a subgraph of the strong product of a graph with bounded treewidth and a path. This result has been the key tool to resolve important open problems regarding queue layouts, nonrepetitive colourings, centered colourings, and adjacency labelling schemes. In this paper, we extend this line of research by utilizing shallow minors to prove analogous product structure theorems for several beyond planar graph classes. The key observation that drives our work is that many beyond planar graphs can be described as a shallow minor of the strong product of a planar graph with a small complete graph. In particular, we show that powers of planar graphs, $k$-planar, $(k,p)$-cluster planar, fan-planar and $k$-fan-bundle planar graphs have such a shallow-minor structure. Using a combination of old and new results, we deduce that these classes have bounded queue-number, bounded nonrepetitive chromatic number, polynomial $p$-centred chromatic numbers, linear strong colouring numbers, and cubic weak colouring numbers. In addition, we show that $k$-gap planar graphs have at least exponential local treewidth and, as a consequence, cannot be described as a subgraph of the strong product of a graph with bounded treewidth and a path.}
		
		\JournalArxiv{}{The planar graph product structure theorem of Dujmovi\'{c}, Joret, Micek, Morin, Ueckerdt, and Wood [J. ACM 2020] states that every planar graph is a subgraph of the strong product of a graph with bounded treewidth and a path. This result has been the key tool to resolve important open problems regarding queue layouts, nonrepetitive colourings, centered colourings, and adjacency labelling schemes. In this paper, we extend this line of research by utilizing shallow minors to prove analogous product structure theorems for several beyond planar graph classes. The key observation that drives our work is that many beyond planar graphs can be described as a shallow minor of the strong product of a planar graph with a small complete graph. In particular, we show that powers of planar graphs, $k$-planar, fan-planar and $k$-fan-bundle planar graphs have such a shallow-minor structure. Using a combination of old and new results, we deduce that these classes have bounded queue-number, bounded nonrepetitive chromatic number, polynomial $p$-centred chromatic numbers, linear strong colouring numbers, and cubic weak colouring numbers. In addition, we show that $k$-gap planar graphs have at least exponential local treewidth and, as a consequence, cannot be described as a subgraph of the strong product of a graph with bounded treewidth and a path.}
	\end{abstract}
	
	
	\section{\large Introduction}
	Many results in structural graph theory describe complex graph classes in terms of simpler classes. For example, Robertson and Seymour's graph minor structure theorem \cite{robertson2003graph} states that every graph in a proper minor-closed class can be described by a tree-decomposition of graphs that are almost embeddable on a surface with bounded Euler genus. Recently, graph products have emerged as a powerful tool in describing graphs in a complex class. Say a graph $H$ is \defn{contained} in a graph $G$ if $H$ is isomorphic to a subgraph of $G$, writen \defn{$H \subsetsim G$}. 	\citet{DJMMUW20} established that every planar graph is contained in the strong product of a graph with bounded treewidth and a path. 

	\begin{thm}[\cite{DJMMUW20}]\label{GPSTintro}
		Every planar graph is contained in $H \boxtimes P$ for some planar graph $H$ of treewidth at most $8$ and for some path $P$.
	\end{thm}

	This breakthrough has been the key tool to resolve several major open problems regarding queue layouts \cite{DJMMUW20}, nonrepetitive colourings \cite{DEJWW2020nonrepetitive}, centred colourings \cite{debski2020improved}, adjacency labelling schemes \cite{bonamy2020shorter,EJM2020universal,DEGJMM2020adjacency}, twin-width \cite{BKW2022bandwidth} and infinite graphs \cite{HMSTW2021universality}. Treewidth measures how similar a graph is to a tree and is an important parameter in algorithmic and structural graph theory (see \cite{bodlaender1993tourist,HW2017tied,reed1997treewidth}). Graphs with bounded treewidth are a simpler class of graphs compared to planar graphs. \Cref{GPSTintro} therefore reduces problems on a complicated class of graphs (planar graphs) to a simpler class of graphs (bounded treewidth). 
	
	Motivated by \Cref{GPSTintro}, \citet{BDJMW2021rowtreewidth} defined the \defn{row treewidth}, $\rtw(G)$, of a graph $G$ to be the minimum treewidth of a graph $H$ such that $G \subsetsim H\boxtimes P$ for some path $P$. \Cref{GPSTintro} says that planar graphs have row treewidth at most $8$. \citet{UWY2021GPST} strengthened \Cref{GPSTintro} by improving the upper bound to $6$. Other graph classes that have bounded row treewidth include graphs with bounded Euler genus~\cite{DJMMUW20,DHHW21}, apex-minor-free graphs~\cite{DJMMUW20}, $(g,d)$-map graphs~\cite{DMW20}, $(g,\delta)$-string graphs~\cite{DMW20}, $(g,k)$-planar graphs~\cite{DMW20}, squaregraphs~\cite{HJMW2022square} and powers of bounded degree planar graphs~\cite{DMW20}.
	
	\JournalArxiv{
	In this paper, we extend this line of research by utilizing shallow minors and graph products to demonstrate that several beyond planar graph classes have bounded row treewidth. The study of beyond planar graphs is a vibrant research topic within the graph drawing community (see the surveys \cite{DLM2019beyond,HT2020beyond}). A shallow minor of a graph is obtained by contracting subgraphs with small radii then deleting vertices and edges. While beyond planar graph classes often contain arbitrarily large complete graphs as minors, they are `sparse' under shallow minors. The key observation that drives our work is that many beyond planar graphs can be described as a shallow minor of the strong product of a planar graph with a small complete graph. In particular, we show that the following beyond planar classes have such a shallow-minor structure: powers of planar graphs (\Cref{SectionClique}); $(g,k)$-planar graphs (\Cref{SectionkPlanar}); $(g,\delta)$-string graphs (\Cref{SectionString}); $(k,p)$-cluster planar graphs (\Cref{SectionkpPlanar});	fan-planar graphs (\Cref{SectionFan}); and $k$-fan-bundle-planar graphs (\Cref{SectionkFanBundle}). To prove these results, we first planarise the given graph. For most of the above graph classes, we use the standard planarisation technique of inserting dummy vertices at crossing points. For fan-planar graphs, however, our planarisation procedure is more substantial and is of independent interest (see \cref{FanShallow}).
	}
	
	\JournalArxiv{}{
	In this paper, we extend this line of research by utilizing shallow minors and graph products to demonstrate that several beyond planar graph classes have bounded row treewidth. The study of beyond planar graphs is a vibrant research topic within the graph drawing community (see the surveys \cite{DLM2019beyond,HT2020beyond}). A shallow minor of a graph is obtained by contracting subgraphs with small radii then deleting vertices and edges. While beyond planar graph classes often contain arbitrarily large complete graphs as minors, they are `sparse' under shallow minors. The key observation that drives our work is that many beyond planar graphs can be described as a shallow minor of the strong product of a planar graph with a small complete graph. In particular, we show that the following beyond planar classes have such a shallow-minor structure: powers of planar graphs (\Cref{SectionClique}); $(g,k)$-planar graphs (\Cref{SectionkPlanar});	$(g,\delta)$-string graphs (\Cref{SectionString}); fan-planar graphs (\Cref{SectionFan}); and $k$-fan-bundle-planar graphs (\Cref{SectionkFanBundle}). To prove these results, we first planarise the given graph. For most of the above graph classes, we use the standard planarisation technique of inserting dummy vertices at crossing points. For fan-planar graphs, however, our planarisation procedure is more substantial and is of independent interest (see \cref{FanShallow}).
	}

	A key contribution of this paper is to show that row treewidth is `well-behaved' under shallow minors (see \Cref{SectionKey}). In particular, we show that if $G$ has bounded row treewidth and $H$ is an $r$-shallow minor of $G \boxtimes K_{\ell}$, then $H$ has bounded row treewidth (see \Cref{GPSTshallow}). We conclude that the above beyond planar graph classes have bounded row treewidth. Note that using the concept of shortcut systems, \citet{DMW20} showed that $(g,k)$-planar graphs, $(g,d)$-map graphs, and $(g,\delta)$-string graphs have bounded row treewidth. However, shortcut systems are limited in that they only apply to graph classes with linear crossing number (see \Cref{ShortcutGapPlanar}). Shallow minors subsume and generalise shortcut systems. In particular, the results for fan-planar graphs use shallow minors in their full generality since fan-planar graphs have super-linear crossing numbers. 
	
	Having established that the above classes have a product structure, we employ tools from graph product structure theory to deduce important structural properties of these classes. In particular, we show that each of these classes has bounded queue-number (\Cref{SectionQueue}), bounded nonrepetitive chromatic number (\Cref{SectionNonrepetitive}), polynomial $p$-centred chromatic numbers (\Cref{SectionCenteredColourings}), linear strong colouring numbers (\Cref{SecionColouringNumbers}), and cubic weak colouring numbers (\Cref{SecionColouringNumbers}).
	
	Motivated by the application of shallow minors in describing graph classes, we further explore the structural properties of shallow minors. Since shallow minors are fundamental to graph sparsity theory (see \cite{nevsetvril2012sparsity}), many papers have considered their properties including the grad function \cite{nevsetvril2012sparsity,HQ2017lowdensity}, separators \cite{DN2016sublinear}, twin-width \cite{BGKTW2021twinsmall}, and layered treewidth \cite{dujmovic2017layered}. We continue this line of research by demonstrating that the queue-number, strong-colouring numbers, and weak-colouring numbers are `well-behaved' under shallow-minors. For example, our result for colouring numbers. which is of independent interest, states that if $H$ is an $r$-shallow minor of $G \boxtimes K_{\ell}$, then
	\begin{equation*}
		\scol_s(H) \leq \ell \scol_{2rs+2r+s}(G) \quad \text{and} \quad
		\wcol_s(H) \leq \ell \wcol_{2rs+2r+s}(G).
	\end{equation*}
	
	Our final results show limitations of the above methods. In particular, we show that $k$-gap planar, fan-crossing free, $3$-quasi planar and right angle crossing (RAC) graphs cannot be described as a shallow minor of the strong product of a planar graph with a small complete graph. Indeed, we show that $1$-gap planar graphs have at least exponential local treewidth and thus have unbounded row treewidth (see \Cref{SectionGap}). For fan-crossing free, $3$-quasi planar and RAC graphs, we draw on known results to conclude that these classes are somewhere dense (see \Cref{SectionSomewhereDense}). As such, they have unbounded row treewidth, unbounded stack-number and unbounded queue-number.

	\section{\large Preliminaries}
	We now introduce the relevant notation and background for this paper. See the textbook by \citet{diestel2017graphtheory} for further background.
	
	\subsection{\normalsize Notation}
	Let $G$ be a finite, simple and undirected graph. 
	For $v \in V(G)$, let $N_G(v):=\{u\in V(G):uv \in E(G)\}$ and $N_G[v]:=N_G(v)\cup \{v\}$. When the graph $G$ is clear, we drop the subscript $G$ and use the notation $N(v)$ and $N[v]$. Let $\Delta(G)$ and $\delta(G)$ respectively denote the maximum and minimum degree of $G$. A \defn{vertex $c$-colouring} of $G$ is any function $\phi: V(G) \to C$ where $|C|\leq c$, and an \defn{edge $c$-colouring} of $G$ is any function $\phi: E(G) \to C$ where $|C|\leq c$. For an integer $k \geq 1$, the \defn{$k$-th power} of $G$, \defn{$G^k$}, is the graph with $V(G^k):=V(G)$ and $uv \in E(G^k)$ if $\dist_G(u,v)\leq k$ and $u \neq v$. The \defn{complement} of $G$ is denoted by $\overline{G}$. Thus $\overline{K_n}$ is the edgeless graph on $n$ vertices. A \defn{class} of graphs $\G$ is a set of graphs that is closed under isomorphism. $\G$ is \defn{hereditary} if it is closed under induced subgraphs and \defn{minor-closed} if it is closed under minors. For integers $a,b$ where $a \leq b$, let $[a,b]:=\{a,a+1,\dots,b\}$ and $[b]:=[1,b]$. 
	
	
	
	\subsection{\normalsize Treewidth}
	For a graph $G$, a \defn{tree-decomposition} of $G$ consists of a collection $\B=\{B_x\sse V(G): x \in V(T)\}$ of \defn{bags} indexed by a tree $T$ such that
	(i) for each $v \in V(G)$, $T[\{x \in V(T): v \in B_x\}]$ is a non-empty subtree of $T$; and 
    (ii) for each $uv \in E(G)$, there is a node $x \in V(T)$ such that $u,v \in B_x$.
	We denote a tree-decomposition by the pair $(T,\B)$. The \defn{width} of $(T,\B)$ is $\max\{|B_x|-1:x \in V(T)\}$. The \defn{treewidth} of a graph $G$, $\tw(G)$, is the minimum width of a tree-decomposition of $G$. Tree-decompositions were introduced by Robertson and Seymour \cite{robertson1986algorithmic}.
	
	
	\subsection{\normalsize Graph Products}	
	Let $G_1$ and $G_2$ be graphs. A \defn{graph product} $G_1 \bullet G_2$ is defined with vertex set $V(G_1 \bullet G_2):=$ $\{(a,v)~:~ a~\in~ V(G_1), v~\in~ V(G_2)\}.$ The \defn{lexicographic product} $G_1~\circ~G_2$ consists of edges of the form $(a,v)(b,u)$ where either $ab \in E(G_1)$, or $a=b$ and $uv \in E(G_2)$. The \defn{strong product} $G_1 \boxtimes G_2$ consists of edges of the form $(a,v)(b,u)$ where either $ab \in E(G_1)$ and $v=u$ or $uv \in E(G_2)$, or $a=b$ and $uv \in E(G_2)$. We frequently make use of the well-known fact that $\tw(G \boxtimes K_n)\leq (\tw(G)+1)n-1$ for every graph $G$ and integer $n\geq 1$ (see \cite{BGHK1995approximating} for an implicit proof). For a graph product $\bullet\in \{\circ, \boxtimes\}$, graph class $\G$ and graph $H$, let $\G \bullet H$ denote the class of graphs isomorphic to graphs in $\{G \bullet H: G \in \G\}$. 
	
	In addition to \Cref{GPSTintro}, \citet{DJMMUW20} proved the following product structure theorem for planar graphs.
	\begin{thm}[\cite{DJMMUW20}]\label{GPST}
		Every planar graph is contained in $H \boxtimes P \boxtimes K_3$ for some planar  graph $H$ of treewidth at most $3$ and for some path $P$.
	\end{thm}
	
	The \defn{Euler genus} of the orientable surface with $h$ handles is $2h$. The \defn{Euler genus} of the non-orientable surface with $c$ cross-caps is $c$. The \defn{Euler genus} of a graph $G$ is the minimum integer $k\geq 0$ such that $G$ embeds in a surface of Euler genus $k$; see \cite{MT2001surfaces} for more about graph embeddings in surfaces. Strengthening an earlier result of \citet{DJMMUW20}, \citet{DHHW21} proved the following product structure theorem for graphs on surfaces.
	\begin{thm}[\cite{DHHW21,DJMMUW20}]\label{GPSTgenus}
		Every graph of Euler genus $g$ is contained in $H \boxtimes P \boxtimes K_{\max\{2g,3\}}$ for some planar graph $H$ of treewidth at most $3$ and for some path $P$.
	\end{thm}

	\subsection{\normalsize Beyond Planar Graphs}
	Beyond planar graphs is a vibrant research topic that studies graph classes defined by drawings that forbid certain crossing configurations. See the recent survey by \citet{DLM2019beyond} as well as the monograph by \citet{HT2020beyond}. A key objective of this paper is to understand the global structure of beyond planar graphs.
	
	An \defn{embedded graph} is a graph $G$ with $V(G)\subset \R^2$ in which each edge $uv \in E(G)$ is a curve in $\R^2$ with endpoints $u$ and $v$ and not containing any vertex of $G$ in its interior. A \defn{crossing} in $G$ is a triple $(p,\{uv,xy\})$ with $p \in \R^2$, and $uv,xy \in E(G)$ such that $p \in (uv \cap xy)\bs \{u,v,x,y\}$. $G$ is \defn{plane} if it has no crossings. $G$ is \defn{simple} if any two edges share at most one point in common, including endpoints. The \defn{crossing number}, $\cross(G)$, of a graph $G$ is the minimum number of crossings in an embedded graph isomorphic to $G$. A graph class $\G$ has \defn{linear crossing number} if there exists a constant $c>0$ such that $\cross(G)\leq c|V(G)|$ for every $G \in \G$. By the Crossing Lemma~\cite{leighton1983VLSI,ACNS1982crossing}, this is equivalent to there being a constant $c'>0$ such that $\cross(G)\leq c'|E(G)|$ for every $G \in \G$. 
	
	We now introduce the beyond planar graphs that are relevant to this paper. An embedded graph $G$ is:
	\begin{compactitem}
		\item \defn{$k$-planar }if each edge of $G$ is involved in at most $k$ crossings \cite{PT1997crossings}; 
		\item \defn{$k$-quasi-planar }if every set of $k$ edges do not mutually cross;
		\item \defn{$k$-gap planar} if every crossing can be charged to one of the two edges involved so that at most $k$ crossings are charged to each edge \cite{BBCEEGHKMRT2018gap}; 
		\item \defn{fan-crossing free} if for each edge $e \in E(G)$, the edges that cross $e$ form a matching \cite{CHKK2015fancrossing}; 
		\item \defn{fan-planar }if for each edge $e\in E(G)$ the edges that cross $e$ have a common end-vertex and they cross $e$ from the same side (when directed away from their common end-vertex) \cite{KU2014density}; or
		\item \defn{right angle crossing }(RAC) if each edge is drawn as a straight line segment and edges cross at right angles \cite{DEL2011RAC}.
	\end{compactitem}
	A graph is \defn{$k$-planar}, \defn{$k$-quasi planar}, \defn{$k$-gap planar}, \defn{fan-crossing free}, \defn{fan-planar}, or \defn{RAC} if it is respectively isomorphic to an embedded graph that is $k$-planar, $k$-quasi planar, $k$-gap planar, fan-crossing free, fan-crossing, or RAC. 
	
	\subsection{\normalsize Shallow Minors and Graph Sparsity}
	Let $G$, $H$ and $J$ be graphs and let $r\geq 0$ be an integer and $s\geq 0$ be a half-integer (that is, $2s$ is an integer). $H$ is a \defn{minor} of $G$ if a graph isomorphic to $H$ can be obtained from $G$ by vertex deletion, edge deletion, and edge contraction. A \defn{model} of $H$ in $G$ is a function $\mu$ with domain $V(H)$ such that (i) $\mu(v)$ is a connected subgraph of $G$; (ii) $\mu(v)\cap \mu(w)=\emptyset$ for all distinct $v,w\in V(G)$; and	(iii) $\mu(v)$ and $\mu(e)$ are adjacent for every edge $vw \in E(H)$.	It is folklore that $H$ is a minor of $G$ if and only if $G$ contains a model of $H$. If there exists a model $\mu$ of $H$ in $G$ such that $\mu(v)$ has radius at most $r$ for all $v \in V(H)$, then $H$ is an \defn{$r$-shallow minor} of $G$. If there is a model $\mu$ of $H$ in $G$ such that $|V(\mu(v))|\leq r$ for all $v \in V(H)$, then $H$ is an \defn{$r$-small minor} of $G$. We say that $J$ is an \defn{$r$-subdivision} of $H$ if $J$ can be obtained from $H$ by replacing each edge $uv$ of $H$ by a path of length at most $r+1$ with end-vertices $u$ and $v$. We say that $H$ is an \defn{$s$-shallow topological minor} of $G$ if a subgraph of $G$ is isomorphic to a $2s$-subdivision of $H$. Note that if $H$ is an $s$-shallow topological minor of a graph $G$, then $H$ is an $r$-shallow minor of $G$ whenever $s\leq r$. 

    
	Shallow minors and shallow topological minors are fundamental to graph sparsity theory (see \cite{nevsetvril2012sparsity}). In particular, they define bounded expansion, nowhere dense, and somewhere dense graph classes. For an integer $r\geq 0$ and graph $G$, let $\nabla_r(G)$ be the greatest average degree of an $r$-shallow minor of $G$. A hereditary graph class $\G$ has \defn{bounded expansion} with \defn{expansion function} $f_{\G}: \N\cup \{0\} \to \mathbb{R}$ if $\nabla_r(G)\leq f_{\G}(r)$ for every $r \geq 0$ and graph $G \in \G$. If $f_{\G}$ is polynomial, then $\G$ has \defn{polynomial expansion}. A hereditary graph class $\G$ is \defn{somewhere dense} if there exists an integer $r\geq 0$ such that every graph $H$ is an $r$-shallow minor of some graph $G \in \G$. If $\G$ is not somewhere dense, then it is \defn{nowhere dense}.
	
	Bounded expansion is a robust measure of sparsity with many characterisations \cite{zhu2009generalized,nevsetvril2008grad,nevsetvril2012sparsity}. For example, a graph class $\G$ has bounded expansion if and only if there exists a function $f$ such that every $r$-shallow minor has chromatic number at most $f(r)$ \cite[Proposition~5.5]{nevsetvril2012sparsity}. Moreover, graph classes with polynomial expansion have been characterised as those that admit strongly sublinear separators \cite{DN2016sublinear}. To prove this, \citet{DN2016sublinear} showed that if a graph $G$ admits strongly sublinear separator, then every $r$-shallow minor (for fixed $r\geq 1$) of $G$ also admits strongly sublinear separators. 
	
	Before proceeding, note that the operation of taking a shallow minor of the product of a graph with a small complete graph has previously been studied within graph sparsity theory. In particular, \citet{HQ2017lowdensity} showed that $\nabla_r(G\circ K_{\ell})\leq 5 \ell^2 (r+1)^2\nabla_r(G)$ for every graph $G$ (improving on an earlier result by \citet{nevsetvril2008grad}).
	
	
	\section{\large Shallow Minors and Graph Products}\label{PropertiesShallowMinorProducts}
	
	This section explores various properties of shallow minors and graph products. For further structural results concerning graph products, see \cite{dvovrak2020notes,HW2021products} as well as the handbook \cite{hammack2011handbook}.
	
	\subsection{\normalsize Shortcut System}
	We begin by examining the relationship between shallow minors and shortcut systems. \citet{DMW20} introduced shortcut systems as a way to prove product structure theorems for various non-minor-closed graph classes. A set $\Path$ of paths in a graph $G$ is a \defn{$(k,d)$-shortcut system} if every path $P \in \Path$ has length at most $k$ and every vertex $v \in V(G)$ is an internal vertex for at most $d$ paths in $\Pcal$. Let \defn{$G^{\Path}$} denote the supergraph of $G$ obtained by adding the edge $uv$ if $\Pcal$ contains a $(u,v)$-path. \citet{DMW20} observed that $k$-planar graphs, $d$-map graphs, and several other classes can be described by applying shortcut systems to planar graphs. Using the following theorem, they deduced that these classes have a product structure.
	
	\begin{thm}[\cite{DMW20}]\label{ShortcutGPST}
		If $G\subsetsim H \boxtimes P \boxtimes K_{\ell}$, for some graph $H$ of treewidth at most $t$ and $\Pcal$ is a $(k,d)$-shortcut system for $G$, then $G^{\Pcal} \subsetsim J \boxtimes P \boxtimes K_{d\ell(k^3+3k)}$ for some graph $J$ of treewidth at most $\binom{k+t}{t}-1$ and some path $P$.
	\end{thm}

	In \Cref{SectionKey}, we adapt the proof of \Cref{ShortcutGPST} to establish an analogous result in the more general setting of shallow minors.
	
	\citet{HW2021treedensities} introduced the following variant of shortcut systems. A set $\Path$ of paths in a graph $G$ is a \defn{$(k,d)^*$-shortcut system} if every path in $\Pcal$ has length at most $k$ and for every $v \in V(G)$, if $M_v$ is the set of vertices $u \in V(G)$ such that there exists a $uw$-shortcut in $\Pcal$ in which $v$ is an internal vertex, then $|M_v|\leq d$. Note that $(k,d)$-shortcut systems and $(k,d)^*$-shortcut systems differ in how they count the number of paths that a vertex contributes to. Every $(k,d)^*$-shortcut system is a $(k,\binom{d}{2})$-shortcut system, and every $(k,d)$-shortcut system is a $(k,2d)^*$-shortcut system. Thus, $(k,d)^*$-shortcut systems can give better bounds compared to $(k,d)$-shortcut systems. Shallow minors inherent this strength of $(k,d)^*$-shortcut systems. 
	
	
	We show that graphs obtained by applying a shortcut system to a planar graph are $k$-gap planar and thus have linear crossing number. Later, we show that fan-planar graphs have super-linear crossing number and thus cannot be described by shortcut systems (see \cref{SectionFan}).
	
	\JournalArxiv{}{
		\begin{lem}[\cite{HW2021shallowArXiv}]\label{ShortcutGapPlanar}
		If $\Path$ is a $(k,d)$-shortcut system of a planar graph $G$, then $G^{\Path}$ is $\big((d-1)(k-1)+2d\big)$-gap planar.
	\end{lem}
    See the arXiv version for details \cite{HW2021shallowArXiv}.}
    \JournalArxiv{
		\begin{lem}\label{ShortcutGapPlanar}
		If $\Path$ is a $(k,d)$-shortcut system of a planar graph $G$, then $G^{\Path}$ is $\big((d-1)(k-1)+2d\big)$-gap planar.
	\end{lem}
	}
\JournalArxiv{
\beginproof\ Assume that $G$ is a plane graph. For some $\epsilon>\delta>0$, for every vertex $x\in V(G)$ and edge $uv\in E(G)$, let 
$$B_x:=\{p\in \R^2:\dist_{\R^2}(p,x)\leq \epsilon\}\quad \text{and}\quad C_{uv}:=\{p\in \R^2:\dist_{\R^2}(p,uv)\leq \delta\}\bs (B_u\cup B_v).$$	
Choosing $\epsilon$ and $\delta$ to be sufficiently small, we may assume that $B_x,B_y,C_{ab}$ and $C_{uv}$ are non-empty and pairwise-disjoint for all $x,y\in V(G)$ and $ab,uv\in E(G)$. For each $xy$-shortcut $P=(x=w_0,w_1,\dots,w_{\ell-1},w_{\ell}=y)\in \Pcal$, embed the edge $xy$ in the region
$$B_{w_0}\cup C_{w_0w_1}\cup B_{w_1}\cup\dots \cup B_{w_{\ell-1}}\cup C_{w_{\ell-1}w_{\ell}}\cup B_{w_{\ell}}$$
so that for all $uv,xy\in E(G^{\Path})$ with corresponding shortcuts $P_1,P_2\in \Pcal$, 
each crossing between $uv$ and $xy$ occurs in some $B_w$ where $w$ is an internal vertex of $P_1$ or $P_2$, and	$uv$ and $xy$ cross at most once in $B_w$ (see \Cref{fig:Shortcuts}).

		\begin{figure}[!htb]
			\centering\includegraphics[width=0.5\textwidth]{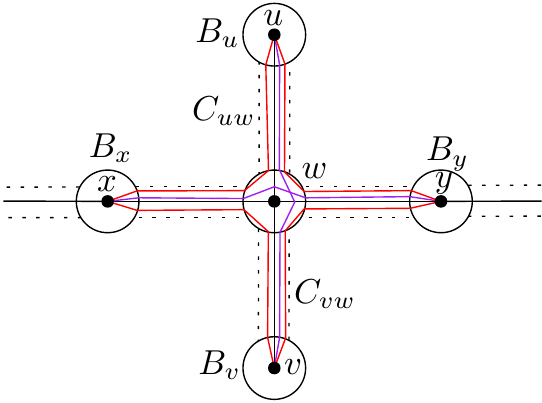}
			\caption{Embedding $G^{\Path}$ into the plane.}
			\label{fig:Shortcuts}
		\end{figure}
		
Let $E_1$ be the edges of $G$ and $E_2:=E(G^{\Pcal})\setminus E_1$ be the new edges in $G^{\Path}$. Since $E_1$ is a set of non-crossing edges, every crossing involves an edge from $E_2$. If $uv\in E_1$ and $xy\in E_2$ cross, then charge the crossing to $uv$. Now suppose $x_1y_1,x_2y_2\in E_2$ cross in $B_w$ for some $w\in V(G)$. Let $P_1,P_2\in \Pcal$ be the shortcuts that respectively correspond to $x_1y_1$ and $x_2y_2$. If $w$ is an internal vertex of $P_1$, then charge the crossing to $x_2y_2$. Otherwise $w$ is an internal vertex of $P_2$, in which case, charge the crossing to $x_1y_1$. 
		
We now upper bound the number of crossings charged to an edge. For an edge $uv \in E_1$, if an edge $xy \in E_2$ crosses $uv$, then the $xy$-shortcut in $\Pcal$ contains $u$ or $v$ as an internal vertex. As such, at most $2d$ crossings are charged to $uv$. Now consider an edge $x_1y_1\in E_2$ with corresponding shortcut $P_1\in \Pcal$. If an edge $x_2y_2\in E_2$ with corresponding shortcut $P_2\in \Pcal$ crosses $x_1y_1$ in $B_w$ and the crossing is charged to $x_1y_1$ for some $w\in V(G)$, then $w$ is an internal vertex of $P_2$ that is contained in $V(P_1)$. Since each vertex is an internal vertex for at most $d$ paths and $V(P_1)$ has at most $(k-1)$ internal vertices, at most $(d-1)(k-1)+2d$ crossings are charged to $x_1y_1$. Therefore, $G^{\Path}$ is $\big((d-1)(k-1)+2d\big)$-gap planar.\qed
}

	We now show that shallow minors generalise shortcut systems.
	
	\begin{lem}\label{ShortcutShallow}
		For every $(k,d)$-shortcut system $\Pcal$ of a graph $G$, $G^{\Path}$ is a $(\frac{k-1}{2})$-shallow topological minor of $G \circ \overline{K_{d+1}}.$
	\end{lem}
	
	\begin{proof}
		Assume that $V(K_{d+1})=[d+1]$ and that for each $uv \in E(G)$ there is a corresponding $(u,v)$-path in $\Pcal$ with length $1$. For each $w \in V(G)$, let $M_w$ be the set of paths in $\Pcal$ that contains $w$ as an internal vertex. Let $\phi_w: M_w \to [2,d+1]$ be an injective function. For each $v\in V(G)$, let $\phi_v(P):=1$ for all shortcuts $P \in \Pcal$ with end-vertex $v$. For each shortcut $P_{uv}=(u=w_0,w_1,\dots,w_{\ell-1},w_\ell=v) \in \Path$ where $\ell \leq k$, let $\tilde{P}_{uv}$ be the path in $G \circ \overline{K_{d+1}}$ defined by $V(\tilde{P}_{uv}):= \{(w_i,\phi_{w_i}(P_{uv})):i \in [0,\ell]\}$ and $E(\tilde{P}_{uv}):= \{(w_i,\phi_{w_i}(P_{uv}))(w_{i+1},\phi_{w_{i+1}}(P_{uv})):i \in [0,\ell-1]\}$. Let $\tilde{\Pcal}$ be the set of such $\tilde{P}_{uv}$ paths.
		
		We claim that $\tilde{\Pcal}$ defines a $(\frac{k-1}{2})$-shallow topological minor of $G^{\Path}$ in $G \circ \overline{K_{d+1}}$ where each vertex $v \in V(G^{\Path})$ is mapped to $(v,1)\in V(G \circ \overline{K_{d+1}})$ and each edge $uv \in E(G^{\Path})$ is mapped to $\tilde{P}_{uv}\in \tilde{\Pcal}$.
		Let $uv \in E(G^{\Path})$. Then there exists a path $P_{uv} \in \Pcal$ with length at most $k$ and end-vertices $u$ and $v$. By construction, $\tilde{P}_{uv}$ is a path in $G \circ \overline{K_{d+1}}$ from $(u,1)$ to $(v,1)$ with length at most $k$. Thus, it suffice to show that the paths in $\tilde{\Pcal}$ are internally disjoint. The internal vertices of each path $\tilde{P}_{uv} \in \tilde{\Pcal}$ are of the form $(w,\phi_w({P}_{uv}))$ where $\phi_w({P}_{uv})\in [2,d+1]$. Suppose there is another path $\tilde{Q}_{xy} \in \tilde{\Pcal}$ for which $(w,\phi_w({P}_{uv}))$ is an internal vertex of $\tilde{Q}_{xy}$. Then $\tilde{P}_{uv}=\tilde{Q}_{xy}$ since $\phi_w$ is injective. As such, the paths in $\tilde{\Pcal}$ are internally disjoint, as required.
	\end{proof}

	\subsection{\normalsize Key Tool}\label{SectionKey}
	Having shown that shallow minors generalise shortcut systems, we now show that shallow minors inherent product structure.
	
	\begin{restatable}{thm}{MainGPST}\label{MainShallowGPST}
		For all graphs $H$ and $L$, if a graph $G$ is an $r$-shallow minor of $H \boxtimes L \boxtimes K_{\ell}$ where $H$ has treewidth at most $t$ and $\Delta(L^r)\leq k$, then $G \subsetsim J \boxtimes L^{2r+1} \boxtimes K_{\ell(k+1)}$ for some graph $J$ with treewidth at most $\binom{2r+1+t}{t}-1$.
	\end{restatable}
	
	To prove \Cref{MainShallowGPST}, we use the language of $H$-partitions from \citet{DJMMUW20}. For a graph $H$, an \defn{$H$-partition} of a graph $G$ is a partition $\mathcal{Y}=(Y_x:x \in V(H))$ of $V(G)$ indexed by the nodes of $H$ such that for any edge $vw \in E(G)$, if $v \in B_x$ and $w \in B_y$, then $xy \in E(H)$ or $x=y$. We say that $H$ is the \defn{quotient} of $\mathcal{Y}$ and is sometimes denoted by $G/\mathcal{Y}$. For our purposes, we require the following extension of $H$-partitions. An \defn{$(H,L)$-partition} $(\mathcal{Y},\mathcal{Z})$ consists of an $H$-partition $\mathcal{Y}$ and an $L$-partition $\mathcal{Z}$ of $G$ for some graphs $H$ and $L$. The \defn{width} of $(\mathcal{Y},\mathcal{Z})$ is $\max\{|Y_y \cap Z_z|: y \in V(H), z \in V(L)\}$. 
	
	The following generalises an observation by \citet{DJMMUW20}.
	
	\begin{obs}\label{OrthogonalPartitions}
		For all graphs $H$ and $L$, a graph $G$ is contained in $H \boxtimes L \boxtimes K_{\ell}$ if and only if $G$ has an $(H,L)$-partition with width at most $\ell$.
	\end{obs}
	
	For a tree $T$ rooted at some node $x_0 \in V(T)$, we say that a node $a \in V(T)$ is a \defn{$T$-ancestor} of $x \in V(T)$ (and $x$ is a \defn{$T$-descendent} of $a$) if $a$ is contained in the path in $T$ from $x_0$ to $x$. If in addition $a\neq x$, then $a$ is a \defn{strict $T$-ancestor} of $x$. Note that $x$ is a $T$-ancestor and $T$-descendent of itself. For each node $x \in V(T)$, define
	$$T_x:=T[\{y \in V(T): \text{$x$ is a $T$-ancestor of $y$}\}]$$
	to be the maximal subtree of $T$ rooted at $x$. 
	
	The proof of \Cref{MainShallowGPST} is an adaptation of the proof of \Cref{ShortcutGPST} \cite[Theorem~9]{DMW20}. We make use of the following well-known normalisation lemma (see \cite[Lemma~2]{DMW20} for a proof).
	
	\begin{lem}\label{NormalisedTdecomp}
		For every graph $H$ of treewidth $t$, there is a rooted tree $T$ with $V(T)=V(H)$ and a width-$t$ tree-decomposition $(T,\{W_x :x \in V(T)\})$ of $H$ that has the following additional properties:
		\begin{compactitem}
			\item[(T1)] for each node $x \in V(H)$, the subtree $T[x]:=T[\{y \in V(T): x \in W_y\}]$ is rooted at $x$; and consequently 
			\item[(T2)] for each edge $xy \in E(H)$, one of $x$ or $y$ is a $T$-ancestor of the other.
		\end{compactitem}
	\end{lem}
	
	
	We now prove our main technical lemma which, together with \Cref{OrthogonalPartitions}, implies \Cref{MainShallowGPST}.	
	\begin{lem}\label{TechnicalLemma}
		Let $G$ be a graph having an $(H,L)$-partition with width $\ell$ in which $H$ has treewidth at most $t$ and $\Delta(L^r)\leq k$. Then every $r$-shallow minor $G'$ of $G$ has a $(J,L^{2r+1})$-partition with width at most $\ell(k+1)$ where the graph $J$ has treewidth at most $\binom{2r+1+t}{t}-1$.
	\end{lem}
	
	\begin{proof}		
		Let $\mu$ be an $r$-shallow model of $G'$ in $G$. Assume that $V(G')\sse V(G)$ and that $u$ is a centre of $\mu(u)$ for each $u \in V(G')$. Let $\mathcal{Y}:=(Y_x:x \in V(H))$ and $\mathcal{Z}:=(Z_x:x \in V(L))$ respectively be an $H$-partition and an $L$-partition of $G$, where $(\mathcal{Y},\mathcal{Z})$ has width at most $\ell$. Let $\Tcal=(T,\B)$ be a normalised tree-decomposition of $H$ that satisfies the conditions of \Cref{NormalisedTdecomp}. For each $x \in V(T)$, let $V_x:=\bigcup_{y \in V(T_x)} Y_y$. Observe that $V_y \sse V_x$ whenever $y$ is a $T$-descendent of $x$. For a vertex $u \in V(G')$, let $X_u:=\{x \in V(T): V(\mu(u)) \cap Y_x \neq \emptyset\}$. 
		
		\textbf{Claim 1:} For every $u \in V(G')$, there exists a node $a(u) \in X_u$ such that $a(u)$ is a $T$-ancestor of every node in $X_u$ and thus, $V(\mu(u)) \sse V_{a(u)}$.
		\begin{proof}
			Since $H$ is a partition of $G$ and $\mu(u)$ is connected, $H[X_u]$ is connected. By the transitivity of the $T$-ancestor relationship, $(T2)$, there exists a node $a(u) \in X_u$ such that $a(u)$ is a $T$-ancestor of every node in $X_u$.
		\end{proof}
		
		For each $x \in V(T)$, define $S_x:=\{u \in V(G'): a(u)=x\}$. Observe that $\mathcal{S}:=(S_x: x \in V(T))$ is a partition of $V(G')$. Let $J:=G'/ \mathcal{S}$ denote the resulting quotient graph, and let $V(J)\sse V(T)$ where each $x \in V(J)$ is obtained by identifying $S_x$ in $G'$. For each $z \in V(L)$, let $Z_z':=Z_z\cap V(G')$ and $\mathcal{Z}':=(Z_z':x \in V(L))$.
		
		From here on in, we need to show: (i) $\mathcal{S}$ has small width with respect to $\mathcal{Z}'$; and (ii), that $J$ has small treewidth. The next claim demonstrates (i).
		
		\textbf{Claim 2:} $|S_x \cap Z_z'|\leq \ell(k+1)$ for all $x \in V(J)$ and $z \in V(L)$.
		\begin{proof}
			Let $u \in S_x \cap Z_z'$ and $W:=\{j \in V(L): \dist_L(z,j)\leq r\}$. By definition, $|W|\leq k+1$. Since $a(u)=x$, there is a vertex $w \in V(\mu(u)) \cap Y_x$. Since $\dist_G(u,w)\leq r$ and $L$ is a partition of $G$, $w \in \bigcup_{j \in W} Z_j'$. As such, $w \in Y_x \cap (\bigcup_{j \in W} Z_{j}')$. Therefore $|S_x \cap Z_z'|\leq \ell(k+1)$ since each vertex in $Y_x \cap (\bigcup_{j \in W} Z_{j}')$ can contribute at most one vertex to $S_x \cap Z_z'$ and $|Y_x \cap (\bigcup_{j \in W} Z_{j}')|\leq \ell(k+1)$.
		\end{proof}
		The following claim will be useful in bounding the treewidth of $J$.
		
		\textbf{Claim 3:} For each edge $xy \in E(J)$, one of $x$ or $y$ is a $T$-ancestor of the other.
		\begin{proof}
			Our goal is to show that there exists a node $w \in V(T)$ such that both $x$ and $y$ are $T$-ancestors of $w$. This implies that one of $x$ or $y$ is a $T$-ancestor of the other.
			
			Since $xy \in E(J)$, $G'$ contains an edge $uv$ with $u \in S_x$ and $v \in S_y$. Let $\tilde{u}\tilde{v}\in E(G)$ where $\tilde{u} \in V(\mu(u))$ and $\tilde{v} \in V(\mu(v))$. 
			Let $\tilde{x} \in V(T)$ be the unique node where $\tilde{u} \in Y_{\tilde{x}}$ and let $\tilde{y} \in V(T)$ be the unique node where $\tilde{v} \in Y_{\tilde{y}}$. Since $\tilde{u} \in V(\mu(u))$ and $\tilde{v} \in V(\mu(v))$, by Claim 1, $x$ is a $T$-ancestor of $\tilde{x}$, and $y$ is a $T$-ancestor of $\tilde{y}$. So if $\tilde{x}=\tilde{y}$, then by setting $w:=\tilde{x}$ we are done. Otherwise, $\tilde{x}\tilde{y}\in E(H)$ since $xy \in E(G)$. By $(T2)$, either $\tilde{x}$ or $\tilde{y}$ is a $T$-ancestor of the other. Without loss of generality, assume that $\tilde{x}$ is a $T$-ancestor of $\tilde{y}$. Then by setting $w:=\tilde{y}$, we are done.
		\end{proof}
		We now show that $J$ has small treewidth.
		
		\textbf{Claim 4:} $J$ has a tree-decomposition in which every bag has size at most $\binom{2r+1+t}{t}$.
		\begin{proof}
		To define the tree-decomposition of $J$, use the same tree $T$ that was used for the tree-decomposition $\Tcal$ of $H$. For each node $x$ of $T$, let $C_x$ be a bag that contains $x$ as well as every $T$-ancestor $a$ of $x$ such that $J$ contains an edge $ax'$ where $x'$ is a $T$-descendent of $x$. Clearly, every vertex is in a bag. Claim 3 ensures that every edge is in a bag. The connectedness of $T[a]:=T[\{x \in V(T):a \in C_x\}]$ follows from the fact that for every edge $x'a \in E(J)$ where $x'$ is a $T$-descendent of $a$, every node $x$ on the $(x',a)$-path in $T$ is a node in $T[a]$. As such, $(C_x: x \in V(T))$ defines a tree-decomposition of $J$. It remains to bound the size of each bag $C_x$.

	\begin{minipage}{0.6\textwidth}
	Consider an arbitrary node $x \in V(T)$ where $x_0,\dots,x_h$ is the path from the root $x_0$ of $T$ to $x_h:=x$. To simplify our notation, let $V_i,Y_i$ and $S_i$ be shorthand for $V_{x_i}, Y_{x_i}$ and $S_{x_i}$ respectively.
	If $x_{\delta} \in C_x$, it is because $x_{\delta}x'\in E(J)$ for some $T$-descendent $x'$ of $x$. This implies that $G'$ contains an edge $uv$ with $u \in S_{x'}$ and $v \in S_{\delta}$. Since $a(v)=x_{\delta}$, there exists $\tilde{v} \in V(\mu(v))\cap Y_{\delta}$. Let $\tilde{u} \in V(\mu(u))$ be such that there exists a path $P=(\tilde{u}=w_0,w_1,\dots,w_p=\tilde{v})$ in $G[V(\mu(v))\cup 
			\{\tilde{u}\}]$ where $p\leq 2r+1$ (see \Cref{fig:BallPath}).
			\end{minipage}
\hfill
	\begin{minipage}{0.4\textwidth}
	\begin{figure}[H]
	\centering\includegraphics[width=0.5\textwidth]{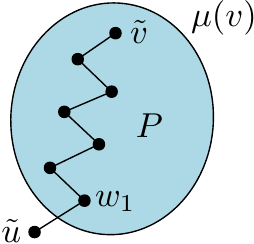}
		\caption{Path $P$ from $\tilde{u}$ to $\tilde{v}$.}
		\label{fig:BallPath}
	\end{figure}
 	\end{minipage}

			For each $i \in [0,p]$, let $s_i:=\max{ \{\ell \in [0,h]:w_0,\dots,w_i \sse V_{\ell}\}}$, and let $a_i:=x_{s_i}$. That is, $a_i$ is the furthest vertex from the root on the $(x_0,x_h)$-path such that $w_0,\dots,w_i \sse V_{a_i}$. Then $s_0,\dots,s_p$ is a non-increasing sequence and $a_0,\dots,a_p$ is a sequence of nodes of $T$ whose distance from the root $x_0$ is non-increasing.
			
			We claim that $a_0=x_h$. Since $a(u)=x_h$, by Claim 1, $\tilde{u} \in V_h$. Since $h$ is trivially the maximum of $\{0,1,\dots,h\}$, $a_0=x_h$ as required.
			
			We claim that $a_p=x_{\delta}$. Since $a(u)$ and $a(v)$ are $T$-descendent of $x_{\delta}$, by Claim 1, $w_0,w_1,\dots,w_p \sse V_{\delta}$ and thus, $a_p\geq x_\delta$. Since $\tilde{v} \in Y_{\delta}$, $a_p \leq x_\delta$ as required. 
			
			Let $H^{+}$ denote the super graph of $H$ with vertex set $V(T)$ in which $xy \in E(H^{+})$ if and only if there exists some $z \in V(T)$ such that $x,y \in B_z$. We claim that $a_0,\dots,a_p$ is a lazy walk in $H^{+}$. Suppose that $a_i\neq a_{i+1}$ for some $i \in [0,p-1]$. Then $w_i \in V_{a_i}$ and $w_{i+1} \not \in V_{a_i}$. Let $b_i$ and $c_i$ be the unique nodes where $w_i \in Y_{b_i}$ and $w_{i+1} \in Y_{c_i}$. By definition of $V_{a_i}$, $b_i$ is a $T$-descendant of $a_i$. Now $b_i c_i \in E(H)$ since $w_i w_{i+1} \in E(G)$, and thus by $(T2)$, $c_i$ is a strict $T$-ancestor of $a_i$. Since $w_1,\dots,w_i,w_{i+1} \in V_{c_i}$ and $w_{i+1} \in Y_{c_i}$, it follows that $a_{i+1}=c_i$. By $(T1)$, $a_{i+1} \in B_{a_{i+1}}$ and $a_{i+1} \in B_{b_i}$. Since $a_i$ is on the path from $b_i$ to $a_{i+1}$ in $T$, this implies $a_{i+1}\in B_{a_i}$. Therefore, $a_ia_{i+1} \in E(H^+)$ as required. 		
		 
			By removing repeated vertices this gives a path of length at most $2r+1$ in the directed graph $\overrightarrow{H}^{+}$ obtained by directing each edge $xy \in E(H^{+})$ from its $T$-descendant $x$ towards its $T$-ancestor $y$. We can then appeal to \cite[Lemma~13]{PS2019centred}, which states that the number of nodes in $\overrightarrow{H}^{+}$ that can be reached from any node $x$ by a directed path of length at most $2r+1$ is at most $\binom{2r+1+t}{t}$. We have shown that for every $x_{\delta} \in C_x$, there is a directed path from $x_{\delta}$ to $x$. It follows that $|C_x|\leq \binom{2r+1+t}{t}$ as required.
		\end{proof}
		To finish the proof, note that $\mathcal{L}$ may not be a valid partition of $G'$. In particular, for every edge $vw\in E(G')$ with $v \in L_i$ and $w \in L_j$, we have $\dist_L(i,j)\in \{0,1,\dots,2r+1\}$. Thus, by indexing $\mathcal{L}$ by $L^{2r+1}$ instead of $L$, we obtain a valid partition of $G'$. Therefore, $G'$ has a $(J,L^{2r+1})$-partition with width at most $(2r+1)(k+1)$ where $\tw(J)\leq \binom{2r+t}{t}-1$.	
	\end{proof}

	Recall that the row treewidth $\rtw(G)$ of a graph $G$ is the minimum treewidth of a graph $H$ such that $G \subsetsim H\boxtimes P$ for some path $P$. Since $\Delta(P^r)\leq 2r$ and $P^{2r+1}\subsetsim P\boxtimes K_{2r+1}$, we have the following consequence of \Cref{MainShallowGPST}:
	
	\begin{thm}\label{GPSTshallow}
		If $G$ is an $r$-shallow minor of $H \boxtimes P \boxtimes K_{\ell}$ where $H$ has treewidth at most $t$ and $P$ is a path, then $G \subsetsim J \boxtimes P \boxtimes K_{\ell(2r+1)^2}$ where $J$ has treewidth at most $\binom{2r+1+t}{t}-1$, and thus $\rtw(G)\leq \binom{2r+1+t}{t}\ell(2r+1)^2-1$.
	\end{thm}
	\cref{GPSTshallow} is used 
	in \cref{SectionShallowBeyond} to establish product structure theorems for beyond planar graphs. 
	

	\subsection{\normalsize Stack and Queue Layouts}\label{SectionQueue}
	Heath, Leighton and Rosenberg~\cite{HLR92,HR1992laying} introduced stack and queue layouts as a way to measure the power of stacks and queues to represent graphs. Let $G$ be a graph and $\preceq$ be a total order on $V(G)$. Two disjoint edges $vw,xy\in E(G)$ with $v\prec w$ and $x\prec y$:
	\begin{compactitem}
		\item \defn{cross} with respect to $\preceq$ if $v\prec x\prec w\prec y$ or $x\prec v\prec y\prec w$;
		\item \defn{nest} with respect to $\preceq$ if $v\prec x\prec y\prec w$ or $x\prec v\prec w\prec y$; and
		\item \defn{overlap} with respect to $\preceq$ if $v=x$ or $w=y$. 
	\end{compactitem}
	Consider a function $\varphi:E(G)\to [k]$ for some integer $k\geq 1$. Then $(\preceq,\varphi)$ is a \defn{$k$-stack layout} of $G$ if $vw$ and $xy$ do not cross for all edges $vw,xy\in E(G)$ with $\varphi(vw) = \varphi(xy)$. Similarly, $(\preceq,\varphi)$ is a \defn{$k$-queue layout} of $G$ if $vw$ and $xy$ do not nest for all edges $vw,xy\in E(G)$ with $\varphi(vw)=\varphi(xy)$. If $vw$ and $xy$ neither nest nor overlap for all edges $vw,xy\in E(G)$ with $\varphi(vw)=\varphi(xy)$, then $(\preceq,\varphi)$ is a \defn{strict $k$-queue layout} of $G$. The minimum integer $s\geq 0$ for which $G$ has an $s$-stack layout is the \defn{stack-number}, $\sn(G)$, of $G$. The minimum integer $q\geq 0$ for which $G$ has a $q$-queue layout is the \defn{queue-number}, $\qn(G)$, of $G$. The minimum integer $q\geq 0$ for which $G$ has a strict $q$-queue layout is the \defn{strict-queue-number}, $\sqn(G)$, of $G$. Stack layouts are equivalent to book embeddings (first defined by \citet{ollmann1973book}) and stack-number is also known as \defn{page-number}, \defn{book-thickness} and \defn{fixed outer-thickness}.

	
	
	\citet{wood2005queue} showed that $\qn(G\boxtimes H)\leq 2\sqn(H)\cdot \qn(G)+\sqn(H)+\qn(G)$ for all graphs $G$ and $H$. Observe that $1 \preceq 2\preceq \dots \preceq \ell$ together with $\phi(ij)=|i-j|$ defines a strict $(\ell-1)$-queue layout of $K_{\ell}$; thus $\sqn(K_{\ell})\leq \ell-1$. Hence, for every graph $G$,
	\begin{equation}
		\qn(G \boxtimes K_{\ell})\leq (2\ell-1) \qn(G)+\ell-1. \tag{$\star$}\label{qnGKl}
	\end{equation}
	\citet{DMW2005layouts} proved that graphs of bounded treewidth have bounded queue-number. The best known bound is $\qn(G)\leq 2^{\tw(G)}-1$, due to \citet{wiechert2017queue}. Using similar techniques to the result of \citet{wood2005queue}, \citet{DJMMUW20} proved the following.
	
	\begin{lem}[\cite{DJMMUW20}]\label{QueueProduct}
		If $G \subsetsim H \boxtimes P \boxtimes K_{\ell}$ then $\qn(G)\leq 3\ell \, 2^{\tw(H)}+\floor{\frac{3}{2}\ell}$.
	\end{lem}
	
	Using \Cref{GPST,QueueProduct}, \citet{DJMMUW20} showed that planar graphs have queue-number at most $49$. This result answered an important open problem of \citet{HLR92} as to whether planar graphs have bounded queue-number. Refining the proof in \citep{DJMMUW20}, \citet{BGR2021queue} strengthened this upper bound to $42$.
	
	For stack layouts, much less is known about their behaviour under products. Recently, it was shown that the stack-number of the strong product of any three sufficiently large connected graphs is arbitrarily large \cite{EHMNSW22}. It is open whether $T\boxtimes P$ has bounded stack-number for every tree $T$ and path $P$ \cite{dujmovic2020stack,pupyrev2020book}. 
	
	Now consider the relationship between shallow minors and stack and queue layouts. For stack layouts, \citet{dujmovic2020stack} described a family of graphs $(G_n)_{n \in \N}$ with unbounded stack-number where for every $n \in \N$, the $5$-subdivision of $G_n$ has stack-number at most $3$. As such, stack-number is not `well-behaved' under shallow topological minors. This result has been recently amplified to show that stack-number is not even `well-behaved' under small minors; that is, there exists a graph class $\G$ with bounded stack-number such that the graph class that consist of $2$-small minors of graphs in $\G$ has unbounded stack-number \cite{EHMNSW22}. On the other hand, queue-number is well-behaved: \citet{DW2005subdivisions} showed that if $H$ is an $r$-shallow topological minor of a graph $G$, then $\qn(H)\leq \frac{1}{2}(2\qn(G)+1)^{4r}-1$. The next lemma generalises this result to the realm of shallow minors.
	
	\begin{lem}\label{QueueShallow}
		For every graph $G$, and for every $r$-shallow minor $H$ of $G$,
		$$\qn(H) \leq 2r(2\qn(G))^{2r}.$$
	\end{lem}
	\begin{proof}
		Let $H$ be an $r$-shallow minor of $G$ where $V(H)\sse V(G)$ and let $\mu$ be an $r$-shallow model of $H$ in $G$ where $v$ is a centre of $\mu(v)$ for each $v\in V(H)$. 
		
		Let $(\preceq_G,\psi_G)$ be a $\qn(G)$-queue layout of $G$. For a directed edge $\overrightarrow{xy}\in E(G)$, we say that it is \defn{forward} if $x \prec_G y$. Let $\preceq_H$ be the induced order of $V(H)$; that is, if $u \prec_G v$ for some $u,v \in V(H)$, then $u \prec_H v$. Let $uv \in E(H)$ where $u \prec_H v$ and $P_{u,v}=(u=w_0,\dots,w_j=v)$ be a $(u,v)$-path in $G$ where $j \leq 2r$ and $P_{uv}\sse \mu(u)\cup \mu(v)$. For each $i \in [j]$, direct the edge $\overrightarrow{w_{i-1}w_i}$. Record the following for the queue assignment $\psi_H(uv)$: the length of $P_{uv}$, $\psi_G(w_{i-1}w_i)$ for each $i \in [j]$, and $\{i \in [j]: \text{$\overrightarrow{w_{i-1}w_i}$ is forwards}\}$. Thus, there is at most $2r \times\qn(G)^{2r}\times2^{2r}$ queue assignments of the edges of $H$.
		
		We claim that $(\preceq_H,\psi_H)$ is a $(2r(2\qn(G))^{2r})$-queue layout of $G$. Suppose there exists $uv,xy\in E(H)$ such that $\psi_H(uv)=\psi_H(xy)$ with corresponding paths $P_{uv}=(u=w_0,\dots,w_j=v)$ and $P_{xy}=(x=z_0,\dots,z_j=y)$ in $G$. Now if $uv$ and $xy$ share a common end-vertex, then they do not nest. So assume that $\{u,v\}\cap \{x,y\}= \emptyset$ and thus $V(P_{uv})\cap V(P_{xy})=\emptyset$. Moreover, assume that $u \prec_H v$ and $x \prec_H y$ and $u \prec_H x$. We claim that $w_i \prec_G z_i$ for all $i \in [0,j]$. For $i=0$, the claim holds by assumption. By induction, we may assume that $w_{i-1}\prec_G z_{i-1}$. Now since $\psi_H(uv)=\psi_H(xy)$, we have $\psi_G(w_{i-1}w_i)=\psi_G(z_{i-1}z_{i})$ and without loss of generality, assume that $\overrightarrow{w_{i-1}w_i}$ and $\overrightarrow{z_{i-1}z_{i}}$ are both forward. Now if $z_{i} \prec_G w_i$, then $w_{i-1} \prec_G z_{i-1} \prec_G z_{i} \prec_G w_i$. Since $\{w_{i-1},w_i\}\cap \{z_{i-1},z_i\}=\emptyset$, it follows that $w_{i-1}w_i$ and $z_{i-1}z_{i}$ nest, a contradiction. Therefore, $v \prec_H y$ and hence $xy$ and $uv$ do not nest, as required.
	\end{proof}
	For the beyond planar graphs defined in \Cref{SectionShallowBeyond}, we can use their product structure (with \cref{QueueProduct}) or their shallow-minor structure (with \Cref{QueueShallow}) to prove that they have bounded queue-number. In some cases, the product structure gives better bounds whereas in other cases, the shallow-minor structure gives better bounds.

	\subsection{\normalsize Nonrepetitive Colourings}\label{SectionNonrepetitive}
	Nonrepetitive colourings were introduced by \citet{AGHR2002nonrepetitive} and have since been widely studied (see the survey \cite{wood2020nonrepetitive}). A vertex $c$-colouring $\phi$ of a graph $G$ is \defn{nonrepetitive} if for every path $v_1,\dots,v_{2h}$ in $G$, there exists $i \in [h]$ such that $\phi(v_i)\neq \phi(v_{i+h})$. The \defn{nonrepetitive chromatic number}, $\pi(G)$, is the minimum integer $c\geq 0$ such that $G$ has a nonrepetitive $c$-colouring. 
	
	\citet{KP2008nonrepetitive} showed that $\pi(G)\leq 4^{\tw(G)}$ for every graph $G$. Building upon this result, \citet{DEJWW2020nonrepetitive} proved the following:
	
	\begin{lem}[\cite{DEJWW2020nonrepetitive}]
	\label{NonrepProduct}
		If $G \subsetsim H \boxtimes P \boxtimes K_{\ell}$ then $\pi(G)\leq \ell \, 4^{\tw(H)+1}$.
	\end{lem}
	
	\citet{DEJWW2020nonrepetitive} used \Cref{GPST,NonrepProduct} to show that planar graphs have bounded nonrepetitive chromatic number, resolving a long-standing conjecture of \citet{AGHR2002nonrepetitive}.
	
	We highlight the following open problem concerning nonrepetitive chromatic number. \citet{NOW2012examples} showed that nonrepetitive chromatic number is `well-behaved' under shallow topological minors. We ask whether this can be generalised to shallow minors. More formally, does there exists a function $f$ such that $\pi(H)\leq f(r,\pi(G))$ whenever $H$ is an $r$-shallow minor of $G$?
	
	\subsection{\normalsize Centred Colourings}
	\label{SectionCenteredColourings}
	
	\citet{nevsetvril2008grad} introduced the following definition. A vertex $c$-colouring $\phi$ of a graph $G$ is \defn{$p$-centred }if, for every connected subgraph $X \sse G$, $|\{\phi(v):v \in V(X)\}|>p$ or there exists some $v \in V(X)$ such that $\phi(v)\neq \phi(w)$ for every $w\in V(X)\bs\{v\}$. That is, $X$ receives more than $p$ colours or some vertex in $X$ receives a unique colour. For an integer $p\geq 1$, the \defn{$p$-centred chromatic number}, $\chi_p(G)$, is the minimum integer $c\geq 0$ such that $G$ has a $p$-centred $c$-colouring. Centred colourings are important within graph sparsity theory as they characterise graph classes with bounded expansion \cite{nevsetvril2008grad}.
	
	
	\citet{debski2020improved} established that $\chi_p(G\boxtimes H)\leq \chi_p(G) \chi(H^p)$ for all graphs $G$ and $H$. \citet[Lemma~15]{PS2019centred} showed that $\chi_p(G)\leq \binom{p+t}{t}$ for every graph $G$ with treewidth at most $t$. Thus we have the following bound for $p$-centred chromatic numbers, first mentioned by \citet{DMW20}.
	\begin{lem}[\cite{debski2020improved,DMW20}]\label{CenteredColouring}
		If $G \subsetsim H \boxtimes P \boxtimes K_{\ell}$ and $\tw(H)\leq t$ then $\chi_p(G)\leq \ell (p+1)\binom{p+t}{t}$.
	\end{lem}

	\subsection{\normalsize Layered Treewidth}
	Layered tree-decompositions were introduced by \citet{dujmovic2017layered} as a precursor to graph product structure theory. A \defn{layered tree-decomposition} $(\mathcal{L},\mathcal{T})$ consists of a layering $\mathcal{L}$ and a tree-decomposition $\mathcal{T}=(T,\B)$ of $G$. The \defn{layered width} of $(\mathcal{L},\mathcal{T})$ is $\max\{|L\cap B|:L \in \mathcal{L},B \in \B\}$. The \defn{layered treewidth}, $\ltw(G)$, of $G$ is the minimum layered width of any layered tree-decomposition of $G$. It is known that $\ltw(G)\leq \rtw(G)+1$ for every graph $G$ \cite{DJMMUW20} but for every integer $k\geq 1$, there exists a graph $G$ with $\ltw(G)= 1$ and $\rtw(G)\geq k$ \cite{BDJMW2021rowtreewidth}. \citet{dujmovic2017layered} showed that every graph $G$ with Euler genus $g$ has $\ltw(G)\leq 2g+3$. In addition, \citet[Lemma~9]{dujmovic2017layered} proved the following:
	
	\begin{lem}[\cite{dujmovic2017layered}]\label{ltwShallow}
		For every graph $G$, and for every $r$-shallow minor $H$ of $G \boxtimes K_{\ell}$,
		$$\ltw(H)\leq \ell(4r+1)\ltw(G).$$
	\end{lem}
	
	We mention several applications of layered treewidth. For strong colouring-numbers (see \Cref{SecionColouringNumbers}), \citet{HW2018improper} proved that $\scol_s(G)\leq \ltw(G)(2s+1)$ for every graph $G$. The \defn{boxicity} of a graph $G$, $\boxicity(G)$, is the minimum integer $d\geq 1$, such that $G$ is the intersection graph of axis-aligned boxes in $\R^d$. \citet{thomassen1986interval} established that planar graphs have boxicity at most $3$. \citet{SW2020boxicity} showed that $\boxicity(G)\leq 6\ltw(G)+4$ for every graph $G$. A closer inspection of their proof gives the following bound. 
	\begin{lem}[\cite{SW2020boxicity}]\label{ltwBoxicity}
		$\boxicity(G)\leq 6\ltw(G)+3$ for every graph $G$.
	\end{lem}
	
	At this stage, it is worth comparing layered treewidth to row treewidth. For many applications, row treewidth has superseded layered treewidth by giving qualitatively stronger bounds. For example, while layered treewidth implies that $n$-vertex planar graphs have queue-number $O(\log (n))$ \cite{dujmovic2017layered}, row treewidth implies that planar graphs have bounded queue-number \cite{DJMMUW20}.\footnote{It is open whether queue-number is bounded by layered treewidth. In fact, \citet{BDJMW2021rowtreewidth} showed that graphs with bounded layered treewidth have bounded queue-number if and only if graphs with layered treewidth $1$ have bounded queue-number.} Nevertheless, there are several applications (in particular, boxicity and strong colouring numbers) where row treewidth has not been shown to give stronger bounds than layered treewidth. Moreover, for many beyond planar graphs, we obtain much better bounds for layered treewidth than we do for row treewidth (see \Cref{SectionShallowBeyond}). For example, using \Cref{ltwShallow,GPSTshallow}, we show that fan-planar graphs have row treewidth at most $1619$ and layered treewidth at most $45$ (see \Cref{SectionFan}). As such, we obtain significantly stronger bounds for the boxicity and strong colouring numbers of fan-planar graphs using layered treewidth than via row treewidth. This highlights the value of layered treewidth, especially for beyond planar graphs.
	
	\subsection{\normalsize Generalised Colouring Numbers}\label{SecionColouringNumbers}
	\citet{KY2003orderings} introduced the following definitions. For a graph $G$, a total order $\preceq$ of $V(G)$, a vertex $v\in V(G)$, and an integer $s\geq 1$, let $R(G,\preceq,v,s)$ be the set of vertices $w\in V(G)$ for which there is a path $v=w_0,w_1,\dots,w_{s'}=w$ of length $s'\in[0,s]$ such that $w\preceq v$ and $v\prec w_i$ for all $i\in[s-1]$, and let $Q(G,\preceq,v,s)$ be the set of vertices $w\in V(G)$ for which there is a path $v=w_0,w_1,\dots,w_{s'}=w$ of length $s'\in[0,s]$ such that $w\preceq v$ and $w\prec w_i$ for all $i\in[s-1]$. For a graph $G$ and integer $s\geq 1$, the \defn{$s$-strong colouring number}, $\scol_s(G)$, is the minimum integer such that there is a total order~$\preceq$ of $V(G)$ with $|R(G,\preceq,v,s)|\leq \scol_s(G)$ for every vertex $v$ of $G$. Likewise, the \defn{$s$-weak colouring number}, $\wcol_s(G)$, is the minimum integer such that there is a total order~$\preceq$ of $V(G)$ with $|Q(G,\preceq,v,s)|\leq \wcol_s(G)$ for every vertex $v$ of $G$. 
	
	Colouring numbers provide upper bounds on several graph parameters of interest. First note that $\scol_1(G)=\wcol_1(G)$ which equals the degeneracy of $G$ plus 1, implying $\chi(G)\leq \scol_1(G)$. A proper graph colouring is \defn{acyclic} if the union of any two colour classes induces a forest; that is, every cycle is assigned at least three colours. The \defn{acyclic chromatic number} $\chi_\text{a}(G)$ of a graph $G$ is the minimum integer $k$ such that $G$ has an acyclic $k$-colouring.
	\citet{KY2003orderings} proved that $\chi_\text{a}(G)\leq \scol_2(G)$ for every graph $G$. Other parameters that can be bounded by strong and weak colouring numbers include game chromatic number \citep{KT1994uncooperative,KY2003orderings}, Ramsey numbers \citep{CS1993ramsey}, oriented chromatic number \citep{KSZ1997acyclic}, arrangeability~\citep{CS1993ramsey}, boxicity \citep{EW2018boxicity}, odd chromatic number \citep{H2022odd} and conflict-free chromatic number \citep{H2022odd}.
	
	Another attractive aspect of strong colouring numbers is that they interpolate between degeneracy and treewidth~\citep{GKRSS2018coveringsnowhere}. As previously noted, $\scol_1(G)$ equals the degeneracy of $G$ plus $1$. At the other extreme, \citet{GKRSS2018coveringsnowhere} showed that $\scol_s(G)\leq \tw(G)+1$ for every integer $s\geq 1$, and indeed 
	$\scol_s(G) \to \tw(G)+1$ as $s\to\infty$.
	
	Colouring numbers are important because they characterise bounded expansion \citep{zhu2009generalized} and nowhere dense classes \citep{GKRSS2018coveringsnowhere}, and have several algorithmic applications \cite{dvorak2014approximation,GKS2017propertiesnowhere}. 
	
	
	We prove the following:
	
	\begin{thm}\label{ShallowCol}
		For every graph $G$, for every $r$-shallow minor $H$ of $G\boxtimes K_{\ell}$, and for every integer $s\geq 1$,
		\begin{equation*}
			\scol_s(H)\leq \ell \scol_{2rs+2r+s}(G) \quad \text{and} \quad	\wcol_s(H) \leq \ell \wcol_{2rs+2r+s}(G).
		\end{equation*}
	\end{thm}
	
	A graph class $\G$ has \defn{linear strong-colouring numbers} if there is a constant $c>0$ such that $\scol_s(G)\leq cs$ for every $G\in\G$ and for every integer $s\geq 1$. Say $\G$ has linear strong colouring numbers. Let $\G_r^{\ell}$ be the class of $r$-shallow minors of graphs in $\G\boxtimes K_{\ell}$. Then for all graphs $H\in\G_r^{\ell}$ and integer $s\geq 1$, 
	$$\scol_s(H) \leq \ell \scol_{2rs+2r+s}(G) \leq c\ell (2rs+2r+s) \leq c\ell (4r+1)s.$$ So $\G_r^{\ell}$ has linear strong colouring numbers (with corresponding constant $c\ell (4r+1)$). In particular, for graphs with bounded Euler genus, van~den~Heuvel~et~al.~\citep{HMQRS2017fixed} proved that $\scol_s(G) \leq (4g+5)s + 2g+1 $ and $\wcol_s(G) \leq \big(2g+\tbinom{s+2}{2}\big)(2s+1)$ for every integer $s\geq 1$ and graph $G$ with Euler genus $g$. Thus \Cref{ShallowCol} implies the following:

	\begin{thm}\label{ColShallowPlanar}
		For every graph $G$ with Euler genus $g$, if $H$ is an $r$-shallow minor of $G \boxtimes K_{\ell}$ then:
		\begin{align*}
			\scol_s(H) &\leq \ell\big( (4g+5)(2rs+2r+s)+2g+1\big)\quad\text{and}\\
			\wcol_s(H) &\leq \ell \big(2g+\tbinom{(2r+1)s+2r+2}{2}\big)((4r+2)s+4r+1).
		\end{align*}	
	\end{thm}
	We use \Cref{ColShallowPlanar} in \Cref{SectionShallowBeyond} to demonstrate that several beyond planar graph classes have linear strong colouring numbers and cubic weak colouring numbers. Layered treewidth also implies linear strong colouring numbers for these graph classes~\cite{HW2018improper}. A strength of our approach is that if the upper bound for the strong or weak colouring numbers for graphs with bounded Euler genus is improved, then this would improve the bounds for the beyond planar graphs in \Cref{SectionShallowBeyond}. For example, \citet{JM2021weak} conjectured that planar graphs have $s$-weak colouring numbers $\Theta(s^2\log(s))$. If true, this would imply that 
	the beyond planar graphs in \Cref{SectionShallowBeyond} have $s$-weak colouring numbers $\Theta(s^2\log(s))$. 
	
	We now work towards proving \Cref{ShallowCol}. Sergey Norin (see \citep{ER2018expansion}) observed that every $r$-shallow minor of a graph $G$ has average degree at most $2 \scol_{4r+1}(G)$. The following lemma generalises this observation (the case $s=1$ implies the above observation):

	\begin{lem}\label{ColShallow}
		For every graph $G$, for every $r$-shallow minor $H$ of $G$, and for every integer $s\geq 1$,
		$$\scol_s(H) \leq \scol_{2rs+2r+s}(G).$$
	\end{lem}
	\begin{proof}
		Let $\preceq$ be a total order of $V(G)$ such that $|R(G,u,\preceq,2rs+2r+s)| \leq \scol_{2rs+2r+s}(G)$ for every vertex $u$ of $G$. 
		Let $\mu$ be an $r$-shallow model of $H$ in $G$. 
		Let $\ell_u$ be the leftmost vertex in $V(\mu(u))$ with respect to $\preceq$. 
		Since the $\mu(u)$'s are pairwise disjoint, the $\ell_u$ are pairwise distinct. 
		Let $\preceq'$ be the total order of $V(H)$ where $u \prec v$ if and only if $\ell_u\prec \ell_v$. 
		Consider a vertex $v$ of $H$ and a vertex $u\in R(H,v,\preceq',s)$. 
		So $u\preceq' v$ and there is a path $(v=w_0,w_1,w_2,\dots,w_t=u)$ of length $t\leq s$ in $H$, such that $v\prec' w_i$ for every $i\in[t-1]$. Thus $\ell_u\preceq \ell_v\prec \ell_{w_i}$ for every $i\in[t-1]$. For $i\in[0,t-1]$ there is an edge $x_iy_{i+1}$ of $G$ where $x_i\in V(\mu(w_i))$ and 
		$y_{i+1}\in V(\mu(w_{i+1}))$. 
		For $i\in[1,t-1]$ there is an $(y_i,x_i)$-path $P_i$ in $\mu(w_i)$ of length at most $2r$ (since $\mu(w_i)$ has radius at most $r$). Similarly, there is an $(\ell_v,x_0)$-path $P_0$ in $\mu(v)$ of length at most $2r$, and there is an $(y_t,\ell_u)$-path $P_{t}$ in $\mu(u)$ of length at most $2r$. 
		Thus $P= \ell_vP_0x_0,y_1P_1x_1,y_2P_2x_2,\dots,y_{t-1}P_{t-1}x_{t-1},y_{t}P_{t}\ell_u$ is an $(\ell_v,\ell_u)$-path in $G$ of length at most $(2r+1)(t+1)-1 \leq (2r+1)(s+1)-1= 2rs+2r+s$. 
		Let $u'$ be the first vertex in $V(\mu(u))$ on $P$ where $u'\preceq \ell_v$. There is such a vertex since $\ell_u$ is a candidate. 
		Let $P'$ be the subpath of $P$ from $\ell_v$ to $u'$. 
		Consider an internal vertex $x$ of $P'$. 
		So $x$ is in $\mu(w_i)$ for some $i\in[t]$.
		If $i=t$ then $\ell_v<x$ by the definition of $u'$. 
		If $i\in[t-1]$ then $\ell_v\prec \ell_{w_i}\preceq x$ by the definition of $\ell_{w_i}$. 
		Thus, every internal vertex of $P'$ is to the right of $\ell_v$ in $\preceq$. 
		So $u'\in R(G,\ell_v,\preceq,2rs+2r+s) \cap V(\mu(u))$. 
		Charge $u$ to $u'$. 
		Each vertex $u\in R(H,v,\preceq',s)$ is charged to a distinct vertex since $u'\in V(\mu(u))$. 
		So $|R(H,v,\preceq',s)| \leq |R(G,\ell_v,\preceq,2rs+2r+s)| \leq \scol_{2rs+2r+s}(G)$ for every vertex $v$ of $H$. 
		Hence $\scol_s(H) \leq \scol_{2rs+2r+s}(G)$.
	\end{proof}
	
	Note that a result of \citet[Corollary~3.5]{zhu2009generalized} implies that there is a function $f$ that bounds the $s$-strong colouring numbers of an $r$-shallow minor $H$ of a graph $G$ as a function of its $f(s)$-strong colouring number. The bound obtained in \Cref{ColShallow} is significantly better.
	
	For weak colouring numbers, \citet{GKRSS2018coveringsnowhere} showed that if $H$ is an $r$-topological minor of $G$ then $\wcol_s(H)\leq 2\wcol_{2rs+s}(G)$ for every integer $s\geq 1$. The following lemma extends this result to shallow minors. We omit the proof since it is analogous to \Cref{ColShallow}.
	
	\begin{lem}
		\label{WeakShallowCol}
		For every graph $G$, for every $r$-shallow minor $H$ of $G$, and for every integer $s\geq 1$, 
		$$\wcol_s(H) \leq \wcol_{2rs+2r+s}(G).$$
	\end{lem}

	For graph products, \citet{dvovrak2020notes} showed that $\scol_s(G\boxtimes H)\leq \scol_s(G)(\Delta(H)+2)^s$ for all graphs $G$ and $H$. Their proof technique also applies to weak-colouring numbers. A closer inspection of their proof reveals the following stronger upper bound.
	
	\begin{lem}[\cite{dvovrak2020notes}]\label{StrongColProduct}
		For all graphs $G$ and $H$ and every integer $s\geq 1$, 
		\begin{align*}
			\scol_s(G\boxtimes H)&\leq \scol_s(G)(\Delta(H^s)+1)\quad\text{and}\\
			\wcol_s(G\boxtimes H)&\leq \wcol_s(G)(\Delta(H^s)+1).
		\end{align*}
	\end{lem}
	\Cref{ShallowCol} follows from \Cref{StrongColProduct,WeakShallowCol,ColShallow} by setting $H=K_{\ell}$. 
	\section{\large Shallow Minors and Beyond Planar Graphs}\label{SectionShallowBeyond}
	This section shows that several beyond planar graph classes can be described as a shallow minor of the strong product of a planar graph with a small complete graph. In fact, we show the slightly stronger result that they are shallow minors of the lexicographic product of a planar graph with a small edgeless graph. We then use the tools from \Cref{PropertiesShallowMinorProducts} to deduce various structural properties of these classes.

\subsection{\normalsize\boldmath Powers of Planar Graphs}\label{SectionClique}
	Recall that the $k$-th power of a graph $G$ is the graph $G^k$ with vertex set $V(G^k):=V(G)$ and $uv \in E(G^k)$ if $\dist_G(u,v)\leq k$ and $u \neq v$. \citet{DMW20} showed that if a graph $G$ has maximum degree $\Delta$, then $G^k=G^{\Path}$ for some $(k,2k\Delta^k)$-shortcut system $\Path$. 
	
	\citet{HW2021treedensities} introduced the following generalisation of low-degree squares of graphs: for a graph $G$ and integer $d \geq 1$, let $G^{(d)}$ be the graph obtained from $G$ by adding a clique on $N_G(v)$ for each vertex $v \in V(G)$ with $\deg_G(v)\leq d$. 	\citet{HW2021treedensities} observed that $G^{(d)}=G^{\Pcal}$ where $\Pcal$ is some $(2,d)^\star$-shortcut system.
	
	We consider a further generalisation of squares of graphs. Let $d \in \N$ and $G$ be a graph. For each vertex $v \in V(G)$, let $M_v\sse N(v)$ where $|M_v|\leq d$ and let $\mathcal{M}:=\{M_v:v\in V(G)\}$. Let $G^{\mathcal{M}}$ denote the graph obtained from $G$ by adding the edges $uw$ to $G$ whenever there exists some $v \in V(G)$ such that $u \in M_v$ and $w \in N(v)$. We call $\mathcal{M}$ a \defn{$d$-clique lift} of $G$. Clearly $G^{(d)}\subsetsim G^{\mathcal{M}}$ for some $d$-clique lift $\mathcal{M}$. 
	 
	\begin{lem}\label{CliqueLift}
		For every graph $G$ with $d$-clique lift $\mathcal{M}$, the graph $G^{\mathcal{M}}$ is a $1$-shallow minor of $G \circ \overline{K_{d+1}}$.
	\end{lem}

	\begin{proof}
		For each $u \in V(G)$, let $S_u:=\{v \in V(G):u \in M_v\}$ and let $\phi_u: M_u \to \{2,\dots,d+1\}$ be an injective function. For each $u\in V(G)$, let $\mu(u)$ be the subgraph of $G\circ \overline{K_{d+1}}$ induced by $\{(u,1)\}\cup \{(w,\phi_v(w)): w \in S_u \}$. We claim that $\mu$ is a $1$-shallow model of $G^{\mathcal{M}}$ in $G\circ \overline{K_{d+1}}.$ 
		
		Let $u,w \in V(G)$ be distinct. First $\mu(u)$ is connected and has radius at most $1$ since $uv \in E(G)$ for all $v \in S_u$. Second, if $S_u\cap S_w=\emptyset$, then $\mu(u)$ and $\mu(w)$ are disjoint. Otherwise, there exists $v \in M_u \cap M_w$. Say $(v,i)\in \mu(u)$ and $(v,j)\in \mu(w)$. Then $i\neq j$ since $\phi_v$ is injective and so $\mu(u)$ and $\mu(w)$ are vertex disjoint.
		
		It remains to show that if $uw \in E(G^{\mathcal{M}})$ then $\mu(u)$ and $\mu(w)$ are adjacent. If $uw \in E(G)$, then $\mu(u)$ and $\mu(w)$ are adjacent since $(u,1)(v,1)\in E(G \circ \overline{K_{d+1}})$. Otherwise, there exists $v\in V(G)$ such that $u,w \in N_G(v)$ and $u \in M_v$ or $w\in M_v$. Assume $u \in M_v$. Then $(v,\phi_v(u))\in V(\mu(u))$ and hence $\mu(u)$ and $\mu(w)$ are adjacent since $(v,\phi_v(u))(w,1)\in E(G \circ \overline{K_{d+1}})$, as required.
	\end{proof}

	Applying \Cref{CliqueLift} with \Cref{GPSTshallow}, we obtain the following product structure theorem for $d$-clique lifts.
	
	\begin{thm}
		For every planar graph $G$, for every $d$-clique lift $\mathcal{M}$ of $G$, the graph $G^{\mathcal{M}} \subsetsim H \boxtimes P \boxtimes K_{27(d+1)}$ for some graph $H$ with treewidth at most $19$, and thus $\rtw(G^{\mathcal{M}})\leq 540(d+1)-1$.
	\end{thm}

	Since clique-lifts generalise squares of graphs, we have the following corollary.
	
	\begin{cor}
		For every planar graph $G$ with maximum degree $\Delta$, the graph
		$G^2 \subsetsim J \boxtimes P \boxtimes
		K_{27(\Delta+1)}$ for some graph $J$ with treewidth at most $19$, and thus $\rtw(G^2) \leq 540(\Delta+1)-1$.
	\end{cor}
	
	More generally, we now show that powers of graphs can be described using shallow minors.
	\begin{lem}
		\label{ShallowMinorPower}
		Let $G$ be a graph and $k\in\mathbb{N}$ and
		$d:=\Delta(G^{\floor{k/2}})$. Then $G^k$ is a
		$\floor{\frac{k}{2}}$-shallow minor of $G \circ \overline{K_{d+1}}$.
	\end{lem}
	
	\begin{proof}
		Let $N_v$ be the set of vertices in $G$ at distance at most
		$\floor{\frac{k}{2}}$ from $v$.
		Say $N_v=\{v_1,\dots,v_{|N_v|}\}$, where $|N_v|\leq d+1$. Then $w\in N_v$ if and only if $v\in N_w$.	For each vertex $v \in V(G)$, let $\mu(v)$ be the subgraph of $G \circ \overline{K_{d+1}}$ induced by $\{(w,i): w \in N_v, v=w_i\}$. It follows that $\mu(v)$ is a connected subgraph of $G\circ \overline{K_{d+1}}$	with radius $\floor{\frac{k}{2}}$, and $\mu(v)\cap \mu(w)=\emptyset$ for distinct $v,w\in V(G)$.
		
		Consider an edge $uv\in E(G^k)$. Let $P_{uv}=(u=w_0,\dots,w_p=v)$ be a shortest $(u,v)$-path in $G$ (and thus $p\leq k$). Then $\dist_G(u,w_{\floor{p/2}})\leq \floor{ \frac{k}{2}}$ and $\dist_G(w_{\floor{p/2}+1},v)\leq \floor{ \frac{k}{2}}$. Thus $(w_{\floor{p/2}},i)\in \mu(u)$ and $(w_{\floor{p/2}+1},j)\in \mu(v)$ for some $i,j\in [d]$, and $(w_{\floor{p/2}},i)(w_{\floor{p/2}+1},j)\in E(G\circ \overline{K_{d+1}})$. Therefore $\mu$ is a $\floor{\frac{k}{2}}$-shallow model of $G^k$ in $G\circ \overline{K_{d+1}}$.
	\end{proof}
	
	
	Applying \Cref{ShallowMinorPower} with	\Cref{MainShallowGPST} by setting $r=\floor{\frac{k}{2}}$, we obtain the following.
	\begin{thm}\label{PowerLGraph}
	Let $G$ be contained in $H \boxtimes L
		\boxtimes K_{\ell}$, where $H$ and $L$ are graphs with $\tw(H)\leq t$,
		and $\ell\in \mathbb{N}$. Let $k \in \N$ and $d:=\Delta(G^{\floor{k/2}})$ and
		$D:=\Delta(L^{\floor{k/2}})$. Then $G^k \subsetsim J \boxtimes L^{2\floor{k/2}+1} \boxtimes K_{\ell(d+1)( D +1)}$ for
		some graph $J$ with treewidth at most $\binom{2\floor{k/2}+t+1}{t}-1$.
	\end{thm}

	Recall that $P^a \subsetsim P \boxtimes K_a$. Thus when $L$ is a path, \Cref{PowerLGraph} implies the following.
		
	\begin{cor}\label{rtwPower}
		Let $G$ be contained in $H \boxtimes P
		\boxtimes K_{\ell}$, where $H$ is a graph with $\tw(H)\leq t$, $P$ is
		a path, and $\ell\in \mathbb{N}$. Let $k\in \N$ and $d:=\Delta(G^{\floor{k/2}})$.
		Then $G^k \subsetsim J \boxtimes P \boxtimes
		K_{\ell(2\floor{k/2}+1)^2(d+1)}$ for some graph $J$ with treewidth at
		most $\binom{2\floor{k/2}+1+t}{t}-1$.
	\end{cor}
	
	Furthermore, \Cref{rtwPower} with \Cref{GPSTshallow} (and $t=\ell=3$) implies:
	
	\begin{cor}\label{PowerProductStructure}
		Let $G$ be a planar graph. Let $k\in\mathbb{N}$ and
		$d:=\Delta(G^{\floor{k/2}})$. Then $G^k \subsetsim J \boxtimes P \boxtimes K_{3(2\floor{k/2}+1)^2(d+1)}$ for some
		graph $J$ with treewidth at most $\binom{2\floor{k/2}+4}{3}-1$.
	\end{cor}
	
	In addition, by applying \Cref{PowerProductStructure} with \Cref{QueueProduct,NonrepProduct,CenteredColouring}, it follows that for every integer $k\geq 1$ and every planar graph $G$, the $k$-power $G^k$ of $G$ has:
	\begin{compactitem}
		\item $\qn(G^k)\leq (9(k+1)^22^{\binom{2\floor{k/2}+4}{3}-1}+\frac{9}{2}(k+1)^2)(\Delta(G^{\floor{k/2}})+1)$;
		\item $\pi(G^k)\leq 3(2\floor{k/2}+1)^2 4^{\binom{2\floor{k/2}+4}{3}}(\Delta(G^{\floor{k/2}})+1)$; and
		\item $\chi_p(G^k)\leq 3(2\floor{k/2}+1)^2(p+1)\binom{p+\binom{2\floor{k/2}+4}{3}-1}{\binom{2\floor{k/2}+4}{3}-1}(\Delta(G^{\floor{k/2}})+1)$.
	\end{compactitem}
	Observe that for all integers $k,\Delta\geq 2$, the complete
	$(\Delta-1)$-ary tree $T$ of height $\lfloor\frac{k}{2}\rfloor$ has
	diameter at most $k$ and maximum degree $\Delta$. Since $T^k$ is a
	complete graph, $\qn(T^k)\geq ~\floor{\frac{|V(T)|}{2}} =\floor{\frac{\Delta(T^{\lfloor{k/2}\rfloor})+1}{2}}$ \cite{HR1992laying} and $\pi(T^k)\geq ~|V(T)| =
	\Delta(T^{\lfloor{k/2}\rfloor})+1$. Therefore, for fixed $k\in
	\mathbb{N}$, the above upper bounds $\qn(G^k)\leqslant~
	O_k(\Delta(G^{\floor{k/2}}))$ and $\pi(G^k)\leqslant~
	O_k(\Delta(G^{\floor{k/2}}))$ on the queue-number and the nonrepetitive chromatic number of $k$-powers of planar graphs $G$ are asymptotically best possible.

	For colouring numbers of $G^k$, \Cref{ShallowMinorPower,ColShallowPlanar} together imply:
	\begin{align*}
		\scol_s(G^k)&\leq (\Delta(G^{\floor{k/2}})+1)(5ks+k+s+1) \quad \text{and}\\ 
		\wcol_s(G^k) &\leq (\Delta(G^{\floor{k/2}})+1)\tbinom{(k+1)s+k+2}{2}((2k+2)s+2k+1).
	\end{align*}
		
	For layered treewidth and boxicity, by applying \Cref{ShallowMinorPower} with \Cref{ltwShallow,ltwBoxicity}, it follows that $\ltw(G^k) \leq 3(\Delta(G^{\floor{k/2}})+1)(2k+1)$ and $\boxicity(G^k)\leq 18(\Delta(G^{\floor{k/2}})+1)(2k+1)+3$.
	
	\subsection{\normalsize\boldmath $(g,k)$-Planar Graphs}\label{SectionkPlanar}
	Recall that a graph $G$ is $k$-planar if it isomorphic to an embedded graph $G'$ in the plane where each edge in $G'$ is involved in at most $k$ crossings.
	This definition has a natural extension for other surfaces. A graph $G'$ embedded on a surface $\Sigma$ is \defn{$(\Sigma, k)$-plane} if every edge of $G'$ is involved in at most $k$ crossings. A graph is \defn{$(g,k)$-planar} if it isomorphic to some $(\Sigma,k)$-plane graph, for some surface $\Sigma$ with Euler genus at most $g$.
	
	\citet{DMW20} observed that every $(g,k)$-planar graph is a subgraph of $G^{\Path}$ for some graph $G$ with Euler genus at most $g$ and some $(k+1,2)$-shortcut system $\Path$. Thus, by \Cref{ShortcutShallow}, every $(g,k)$-planar graph is a $\frac{k}{2}$-shallow-topological minor of $G \circ \overline{K_3}$ where $G$ has Euler genus at most $g$. We obtain a slightly stronger bound by the standard planarisation method. 
	
    \JournalArxiv{}{
    		\begin{lem}[\cite{HW2021shallowArXiv}]\label{kplanarShallow}
        	Every $(g,k)$-planar graph $G$ is a $\frac{k}{2}$-shallow-topological minor of $H \circ \overline{K_2}$ where $H$ has Euler genus at most $g$.
    	\end{lem}
    See the arXiv version for details \cite{HW2021shallowArXiv}.}

	\JournalArxiv{\begin{lem}\label{kplanarShallow}
		Every $(g,k)$-planar graph $G$ is a $\frac{k}{2}$-shallow-topological minor of $H \circ \overline{K_2}$ where $H$ has Euler genus at most $g$.
	\end{lem} 

	\begin{proof}
		Embed $G$ into a surface $\Sigma$ with Euler genus at most $g$ such that every edge of $G$ is involved in at most $k$ crossings. For each crossing $(p,\{uv,xy\})$ where $uv,xy \in E(G)$, insert a dummy vertex $w$ at $p$ to obtain a graph $H$ with Euler genus at most $g$. Let $W:=V(H)\bs V(G)$ be the set of dummy vertices. Each $uv \in E(G)$ corresponds to a path $P_{uv}$ with length at most $k$ in $H$ where each internal vertex of $P$ is a dummy vertex. Let $\Pcal$ be the set of such paths. For each $w \in W$, let $M_w$ be the set of paths in $\Pcal$ that contains $w$ as an internal vertex. Then $|M_w|\leq 2$. Let $\phi_w:M_w\to \{1,2\}$ be an injective function. For each $v \in V(G)$, let $\phi_v(P):=1$ for all paths $P_{uv}$ with end-vertex $v$. For each path $P_{uv}=(u=w_0,w_1,\dots,w_{j-1},w_j=v) \in \Path$ where $j \leq k+1$, let $\tilde{P}_{uv}$ be the path in $H \circ \overline{K_2}$ defined by $V(\tilde{P}_{uv}):=\{(w_i,\phi_{w_i}(P_{uv})):i \in [0,j]\}$ and $E(\tilde{P}_{uv}):=\{(w_i,\phi_{w_i}(P_{uv}))(w_{i+1},\phi_{w_{i+1}}(P_{uv})):i \in [0,j-1]\}$. Let $\tilde{\Pcal}$ be the set of such $\tilde{P}_{uv}$ paths.
		
		We claim that $\tilde{\Pcal}$ defines a $\frac{k}{2}$-shallow topological minor of $G$ in $H \circ \overline{K_2}$ where each vertex $v \in V(G)$ is mapped to $(v,1)\in V(H \circ \overline{K_2})$ and each edge $uv \in E(G)$ is mapped to $\tilde{P}_{uv}\in \tilde{\Pcal}$.
		
		Let $uv \in E(G)$. Then there is a path $P_{uv} \in \Pcal$ with length at most $k+1$ and end-vertices $u$ and $v$. By construction, $\tilde{P}_{uv}$ is a path in $H \circ \overline{K_2}$ from $(u,1)$ to $(v,1)$ with length at most $k$. Thus, it suffice to show that the paths in $\tilde{\Pcal}$ are internally disjoint. Now for a path $\tilde{P}_{uv} \in \tilde{\Pcal}$, its internal vertices are of the form $(w,\phi_w({P}_{uv}))$ where $w \in W$. Now if there is another path $\tilde{Q}_{xy} \in \tilde{\Pcal}$ for which $(w,\phi_w({P}_{uv}))$ is an internal vertex of $\tilde{Q}_{xy}$, then $\tilde{P}_{uv}=\tilde{Q}_{xy}$ since $\phi_w$ is injective. As such, the paths in $\tilde{\Pcal}$ are internally disjoint, as required.
	\end{proof}
	}
	
	Applying \Cref{kplanarShallow} with \Cref{GPSTshallow}, we obtain the following product structure theorems for $(g,k)$-planar graphs. 
	
	\begin{thm}\label{gkplanar}
		Every $(g,k)$-planar graph $G$ is contained in $H\boxtimes P\boxtimes K_{2\max\{2g,3\}(k+1)^2}$ for some graph $H$ with treewidth at most $\binom{k+4}{3}-1$, and thus $G$ has row treewidth at most $2\max\{2g,3\}(k+1)^2\binom{k+4}{3}-1$.
	\end{thm}
	Note that \citet{DMW20} proved that every $(g,k)$-planar graph is a subgraph of $H \boxtimes P \boxtimes K_{\max\{2g,3\}(6k^2+16k+10)}$, for some graph $H$ of treewidth at most $\binom{k+4}{3}-1$. Thus, our results only improve those of \citet{DMW20} by a constant factor. We do not state the bounds that we obtain for the queue-number, nonrepetitive chromatic number and $p$-centred chromatic number as the improvement is not substantial.
	
%
	
	For colouring numbers, applying \Cref{kplanarShallow} with \Cref{ColShallowPlanar}, it follows that for every $k$-planar graph $G$:
	\begin{align*}
		\scol_s(G) &\leq 10(k+1)s+20k+2\quad\text{and}\\
		\wcol_s(G) &\leq \tbinom{(k+1)s+k+2}{2}(4(k+1)s+4k+2),
	\end{align*}
	and for every $(g,k)$-planar graph $G$:
	\begin{align*}
		\scol_s(G) &\leq (8g+10)(k+1)s+(8g+10)k+4g+2\quad\text{and}\\
		\wcol_s(G) &\leq \big(2g+\tbinom{(k+1)s+k+2}{2}\big)((2k+4)s+4k+2).
	\end{align*}
	
	These results for strong colouring numbers differ from the previous bounds by \citet{HW2018improper} (obtained via layered treewidth) by a small constant factor. More importantly, these are the first known polynomial bound on weak colouring numbers for $k$-planar and $(g,k)$-planar.

	\subsection{\normalsize String Graphs}\label{SectionString}
	A \defn{string graph} is the intersection graph of a set of curves in the plane with no three curves meeting at a single point. For an integer $\delta \geq 2$, if each curve is involved in at most $\delta$ intersections with other curves, then the corresponding string graph is called a \defn{$\delta$-string graph}. A \defn{$(g,\delta)$-string graph} is defined analogously for curves on a surface with Euler genus at most $g$.
	
	\citet{DMW20} showed that every $(g,\delta)$-string graph $G$ is a subgraph of $H^{\Path}$ for some graph $G$ with Euler genus at most $g$ and some $(\delta+1,\delta+1)$-shortcut system $\Path$. By \Cref{ShortcutShallow}, $G$ is a $\frac{\delta}{2}$-shallow topological minor of $H\circ \overline{K}_{\delta+2}$. We obtain a stronger bound by the standard planarisation method.
	
	\JournalArxiv{}{
	\begin{lem}	[\cite{HW2021shallowArXiv}]\label{StringShallow}
		Every $(g,\delta)$-string graph $G$ is a $\floor{\frac{\delta}{2}}$-shallow minor of $H \circ \overline{K_2}$ for some graph $H$ with Euler genus at most $g$.
	\end{lem}
    See the ArXiv version for details \cite{HW2021shallowArXiv}.}
	
	\JournalArxiv{
	\begin{lem}\label{StringShallow}
		Every $(g,\delta)$-string graph $G$ is a $\floor{\frac{\delta}{2}}$-shallow minor of $H \circ \overline{K_2}$ for some graph $H$ with Euler genus at most $g$.
	\end{lem}
	
	\begin{proof}
		Let $\mathcal{C}:=\{C_v:v\in V(G)\}$ be a set of curves in a surface of Euler genus at most $g$ such that each curve is involved in at most $\delta$ intersections with other curves and the intersection graph is isomorphic to $G$. For each pair of curves in $\mathcal{C}$ that intersects, add a vertex at their intersection point. Add an edge between each pair of consecutive vertices along a curve in $\mathcal{C}$. If a curve $C_v$ is involve in one crossing, add another vertex $c_v$ to the curve adjacent to the intersection point that $C_v$ contains. If a curve $C_v$ is not involved in any crossing, add a vertex $c_v$ to the curve. Let $H$ be the planar graph obtained and for each vertex $v \in V(G)$, let $M_v$ be the set of vertices in $H$ that are on the curve $C_v$. 
		
		Observe that $H[M_v]$ has radius at most $\floor{\frac{\delta}{2}}$ for each vertex $v \in V(G)$. Moreover, each vertex $w\in V(H)$ is contained in at most two curves. For each $w \in V(H)$, let $\phi_w: \{u,v\} \to \{1,2\}$ be an injective function where $u$ and $v$ index the curves that contains $w$. For each vertex $v \in V(G)$, let $\mu(v)$ be the subgraph of $H \circ \overline{K_2}$ induced by $\{(w,\phi_w(u)): w \in M_v \}$.
		
		We claim that $\mu$ defines a $\floor{\frac{\delta}{2}}$-shallow model of $G$ in $H \circ \overline{K_2}.$ 
		
		Let $u,v \in V(G)$ be distinct. First, since $H[M_v]$ has radius at most $\floor{\frac{\delta}{2}}$ and is connected, so is $\mu(u)$. Second, if $u$ and $v$ do not intersect, then $\mu(u)$ and $\mu(v)$ are disjoint. Otherwise, there exists $w \in C_u \cap C_v$. In which case, $\phi_w(u)\neq \phi_w(v)$ since $\phi_w$ is injective and hence $\mu(u)$ and $\mu(v)$ are vertex disjoint. Finally, if $uv \in E(G)$ then there is a vertex $w_1 \in M_u\cap M_v$. Let $w_2\in M_v$ be a neighbour of $w_1$. Then $(w_1,\phi_{w_1}(u))(w_2,\phi_{w_2}(v))\in E(H \circ \overline{K_2})$ so $\mu(u)$ and $\mu(v)$ are adjacent. Therefore $\mu$ defines a $\floor{\frac{\delta}{2}}$-shallow model of $G$ in $H \circ \overline{K_2}.$
	\end{proof}
	}
	
	Applying \Cref{StringShallow} in conjunction with \Cref{GPSTshallow}, we deduce the following product structure theorem for string graphs.
	
	\begin{thm}\label{gdStringGPST}
		Every $(g,\delta)$-string graph $G$ is contained in $ H\boxtimes P\boxtimes K_{2\max\{2g,3\}(\delta+1)^2}$ for some graph $H$ with treewidth at most $\binom{2 \floor{\delta/2}+4}{3}-1$, and thus $G$ has row treewidth at most $2\max\{2g,3\}(\delta+1)^2\binom{2 \floor{\delta/2}+4}{3}-1$.
	\end{thm}

	Note that \citet{DMW20} previously showed that every $(g,\delta)$-string graph is contained in $H\boxtimes P\boxtimes K_{\max\{2g,3\}(\delta^4+4\delta^3+9\delta^2+10\delta+4)}$ for some graph $H$ with treewidth at most $\binom{\delta+4}{3}-1$. Thus, \Cref{gdStringGPST} improves this result by a factor of $\delta^2$ in the $K_{\ell}$ term. 
	
	For colouring numbers, by applying \Cref{StringShallow} with \Cref{ColShallowPlanar}, it follows that for every $(g,\delta)$-string graph $G$:	
	\begin{align*}
		\scol_s(G) &\leq (8g+10)(\delta+1)s+(8g+10)\delta+4g+2\quad\text{and}\\
		\wcol_s(G) &\leq 2 \big(2g+\tbinom{(\delta+1)s+\delta+2}{2}\big)((2\delta+2)s+2\delta+1).
	\end{align*}

	%
	%
	
	\JournalArxiv{
	\subsection{\normalsize\boldmath $(k,p)$-Cluster Planar Graphs}\label{SectionkpPlanar}

	\citet{DLLRT2019hybrid} introduced the following class of beyond planar graphs. A graph $G$ is \defn{$(k,p)$-cluster planar}\footnote{Note that \citet{DLLRT2019hybrid} called $G$ a \defn{$(k,p)$-planar graph}. We use the language of $(k,p)$-cluster planar graphs to avoid ambiguity with $(g,k)$-planar graphs.} if $V(G)$ can be partitioned into clusters $C_1,\dots,C_n$ where $|C_i|\leq k$ for all $i \in [n]$ such that $G$ admits a drawing, called a \defn{$(k,p)$-cluster planar representation}, where:
	\begin{compactenum}
		\item each cluster $C_i$ is associated with a closed, bounded planar region $R_i$, called a \defn{cluster region};
		\item cluster regions are pairwise disjoint;
		\item each vertex $v \in V(G)$ is identified with at most $p$ distinct points, called \defn{ports}, on the boundary of its cluster region;
		\item each inter-cluster edge $uv\in E(G)$ is a curve that joins a port of $u$ to a port of $v$; and
		\item inter-cluster edges do not cross or intersect cluster regions except at their endpoints.
	\end{compactenum}
	We now show that $(k,p)$-cluster planar graphs have a simple product structure with no dependency on $p$.
	\begin{lem}\label{kpPlanarSubgraph}
		Every $(k,p)$-cluster planar graph $G$ is contained in $ H \boxtimes K_k$ for some planar graph $H$.
	\end{lem}

\begin{proof}
Begin with a $(k,p)$-cluster planar representation of $G$. Replace each cluster region $R_i$ by a vertex $x_i$ such that $x_i$ and $x_j$ are adjacent if there is an edge joining a port on the boundary of $R_i$ to a port on the boundary of $R_j$. In doing so, we obtain a plane graph $H$ (see \Cref{fig:kpPlanar}). For each cluster $C_i$, let $\phi_i:C_i \to [k]$ be an injective function. For each cluster $C_i$ and vertex $v\in C_i$ let $\phi(v)=(x_i,\phi_i(v))\in V(H\boxtimes K_k)$. Observe that $\phi$ is an injective function whose domain is $V(G)$.
				 
\begin{figure}[!htb]
	\centering	\includegraphics[width=0.4\linewidth]{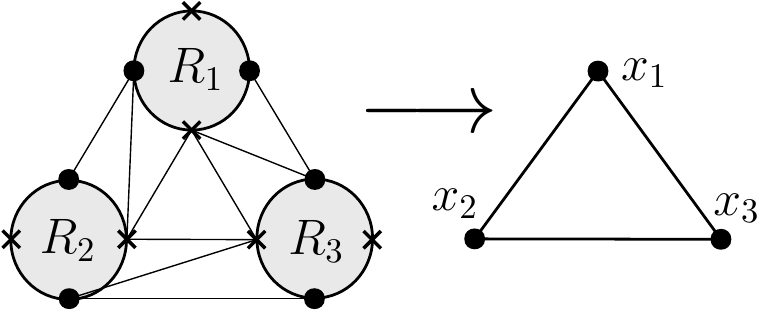}
	\caption{Planarising a $(k,p)$-cluster planar graph.}
	\label{fig:kpPlanar}
\end{figure}
		
We claim that $\phi(u)\phi(v)\in E(H \boxtimes K_k)$ whenever $uv \in E(G)$. First, if $u,v\in C_i$ then $\phi(u)\phi(v)=(x_i,\phi_i(u))(x_i,\phi_i(v))\in E(H \boxtimes K_k)$. Otherwise, if $u \in C_i$ and $v \in C_j$ where $i\neq j$, then $x_ix_j\in E(H)$ since there is an inter-cluster edge from $C_i$ to $C_j$. As such, $\phi(u)\phi(v)=(x_i,\phi_i(u))(x_j,\phi_j(v))\in E(H \boxtimes K_k)$. Thus $G \subsetsim H \boxtimes K_k$, as required. 
	\end{proof}

	Since $(k,p)$-cluster planar graphs can be described by a product, the results from \Cref{PropertiesShallowMinorProducts} apply to this class of graphs. We omit the details.

	We highlight the following consequence of \Cref{kpPlanarSubgraph}. An embedded graph is \defn{IC-planar }if every edge is involved in at most one crossing and the set of all edges that cross form a matching. A graph $G$ is \defn{IC-planar} if it is isomorphic to an IC-planar embedded graph. \citet{DLLRT2019hybrid} observed that the class of IC-planar graphs correspond to $(4,1)$-cluster planar graphs. Thus, by \Cref{kpPlanarSubgraph}, every IC-planar graph $G$ is a subgraph of $H\boxtimes K_4$ for some planar graph $H$.

	%
}	

\subsection{\normalsize Fan-Planar Graphs}\label{SectionFan}

\begin{figure}[!htb]
	\begin{center}
		\includegraphics[width=0.38\linewidth]{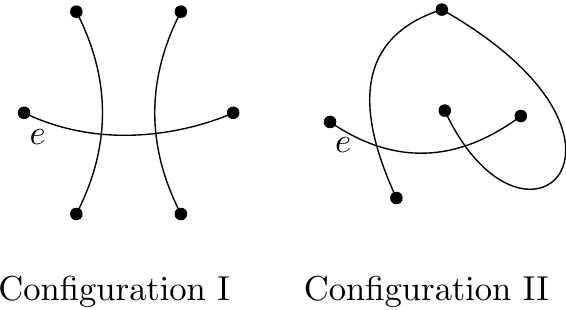}
		\caption{Forbidden crossing configurations for fan-planar graphs.}
		\label{fig:fancrossingforbidden}
	\end{center}
	\vspace*{-0.3cm}
\end{figure}

Recall that an embedded graph $G$ is {fan-planar} if for each edge $e\in E(G)$, the edges that cross $e$ have a common end-vertex and they cross $e$ from the same side (when directed away from their common end-vertex). Equivalently, an embedded graph is fan-planar if it forbids the two crossing configurations in \cref{fig:fancrossingforbidden}. Fan-planar graphs were introduced by \citet{KU2014density}. The class of fan-planar graphs extends $1$-planar graphs and is a proper subset of $3$-quasi-planar graphs.


For an embedded fan-planar graph $G$ and a crossed edge $uv \in E(G)$, a common end-vertex of the edges that cross $uv$ is called a \defn{friend} of $uv$. If $uv$ is crossed at least twice, then it has one friend; otherwise, $uv$ is crossed once and has two friends. A \defn{friend assignment} of $G$ assigns a friend to each crossed edge. For a given friend assignment, a crossed edge $uv$ is \defn{well-behaved} if there exists a non-crossing point $p\in uv$ such that $u$ is the assigned friend of each edge that crosses $uv$ between $p$ and $u$, and $v$ is the assigned friend of each edge that crosses $uv$ between $p$ and $v$. If every crossed edge in $G$ is well-behaved, then the friend assignment is \defn{well-behaved}.

\begin{lem}\label{FriendsFan}
	Every simple fan-planar graph $G$ has a well-behaved friend assignment.
\end{lem}

\beginproof\ 	 Let $uv\in E(G)$ be a crossed edge and let $e_1,\dots,e_m$ be the edges that cross $uv$ ordered by their crossing point from $u$ to $v$. Let $w$ be a common end-vertex of $e_1,\dots,e_m$. If none of $e_1,\dots,e_m$ cross another edge incident to $v$, then let $u$ be the assigned friend of $e_1,\dots,e_m$ and choose $p$ to be an arbitrary point (along $uv$) between $e_m\cap uv$ and $v$. Otherwise, let $i \in [m]$ be minimum such that $e_i$ crosses another edge $vz$ incident to $v$. Choose $p$ to be an arbitrary point on $uv$ between $uv\cap e_{i-1}$ and $uv\cap e_i$. Since none of the edges $e_1,\dots,e_{i-1}$ cross another edge incident to $v$, we may let $u$ be their assigned friend.
	
	
\begin{wrapfigure}{R}{0.5\textwidth}
\hfill
	\includegraphics[width=0.8\linewidth]{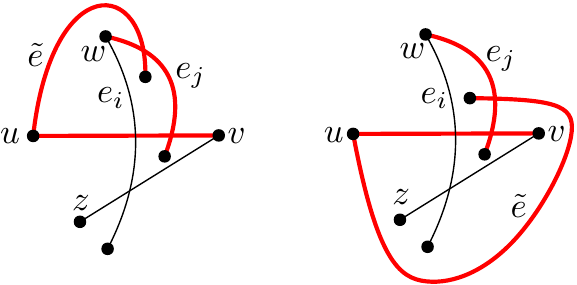}
	\caption{Proof of \cref{FriendsFan}.}
	\label{fig:FanCrossingSide}
\end{wrapfigure}
It remains to show that for each $j \in [i,m]$, $v$ is a common end-vertex of the edges that cross $e_j$ and thus we may let $v$ be the assigned friend of $e_j$. If $e_j$ is only crossed by $uv$, then we are done. If $e_j$ crosses $vz$ then $v$ is $e_j$'s only friend. Otherwise, the end-vertex of $e_j$ 
	opposite to $w$ is contained in the region bounded by the edges $uv,vz,e_i$ (since $G$ is simple fan-planar). Now suppose that another edge $\tilde{e}$ incident to $u$ crosses $e_j$. Since $\tilde{e}$ does not cross $e_i$ or $vz$ (as Configuration I is forbidden), $\tilde{e}$ must cross $e_j$ from the opposite side that $uv$ cross $e_j$, contradicting $G$ being fan-planar (see \cref{fig:FanCrossingSide}). Thus $uv$ is the only edge incident to $u$ that crosses $e_j$ and so $v$ is a common end-vertex of the edges that cross $e_j$, as required.\qed

We now prove our main technical lemma of this subsection.

\begin{lem}\label{FanShallow}
	Every fan-planar graph $G$ is a $1$-shallow minor of $H \circ \overline{K_{3}}$ for some planar graph $H$.
\end{lem} 

\begin{proof}
	By \cite[Theorem~1]{KKRS2021fan} and \cref{FriendsFan}, we may assume that $G$ has a simple fan-planar embedding with a well-behaved friend assignment. Initialise $G^{(0)}:=G$. Arbitrarily choose an edge $u_1v_1\in E(G^{(0)})$ that is involved in at least two crossings and let $w_1 \in V(G^{(0)})$ be the assigned friend of $u_1v_1$. If no such edge exists, then $G$ is $1$-planar and we are done by \cref{kplanarShallow}. Let $E_{1}$ be a maximal subset of $E(G^{(0)})$ such that:
	\begin{compactitem}
		\item $u_1v_1 \in E_{1}$;
		\item each edge in $E_{1}$ has $u_1$, $v_1$ or $w_1$ as its assigned friend; and
		\item $E_1 \setminus V(G^{(0)})$ is a connected subset of the plane.
	\end{compactitem}
	
	Observe that every edge in $E_{1}$ crosses another edge in $E_{1}$ and thus every edge in $E_{1}$ is incident to $u_1,v_1$ or $w_1$.
	
	Let $V_{1}$ be the set of vertices incident to edges in $E_{1}$. Consider the subgraph of $G^{(0)}$ with vertex set $V_1$ and edge set $E_1$. Add a dummy vertex at each crossing point to planarise the subgraph. Let $\tilde{G}_1$ be the plane graph obtained and let $\tilde{D}_1$ be the set of dummy vertices added. Let $G_1':= \tilde{G}_1\cup (G^{(0)}- E_1)$; see \cref{fig:FanDummy}.
	
	\begin{figure}[!htb]
		\begin{center}
			\includegraphics[width=0.3\linewidth]{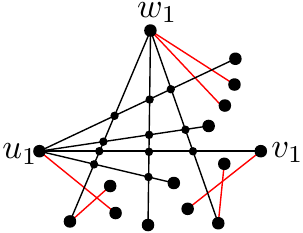}
			\caption{The graph $G_1':= \tilde{G}_1\cup (G^{(0)}- E_1)$ where black edges are from $E(\tilde{G}_1)$ and red edges are from $E(G^{(0)})- E_1$.}
			\label{fig:FanDummy}
		\end{center}
	\end{figure}
	
	\textbf{Claim 1:} If an edge $e\in E(\tilde{G}_1)$ is crossed in $G_1'$, then $e$ is incident to some $z\in V_1\bs \{u_1,v_1,w_1\}$.
	\begin{proof}
		Since $\tilde{G}_1$ is plane, $e$ is crossed by some edge $e'\in E(G^{(0)})\setminus E_1$. Now $e$ is a segment of some edge $xz\in E_1$ where $x \in \{u_1,v_1,w_1\}$ and $z \in V_1$. For the sake of contradiction, suppose that $e$ is not incident to any vertex in $V_1\bs \{u_1,v_1,w_1\}$. Then there exists an edge $\tilde{e}\in E_1$ and a dummy vertex $d=(xz\cap \tilde{e}) \in \tilde{D}_1$ such that $e$ is between $x$ and $d$ (along $xz$). Now if the assigned friend of $\tilde{e}$ is $x$, then the assigned friend of $e'$ is also $x$ since $xz$ is well-behaved. Otherwise, if the assigned friend of $\tilde{e}$ is $z$ then $z\in \{u_1,v_1,w_1\}\setminus\{x\}$. Thus, the assigned friend of $e'$ is in $\{x,z\}\sse \{u_1,v_1,w_1\}$. However, this contradicts the maximality of $E_1$ since $e'$ crosses an edge in $E_1$ and the assigned friend of $e'$ is in $\{u_1,v_1,w_1\}$, as required.
	\end{proof}
	
	For $\epsilon>\delta>0$, for every vertex $x\in V(G_1')$ and edge $uv\in E(G_1')$, let 
	$$B_x:=\{p\in \R^2:\dist_{\R^2}(p,x)\leq \epsilon\}\quad \text{and}\quad C_{uv}:=\{p\in \R^2:\dist_{\R^2}(p,uv)\leq\delta\}\bs (B_u\cup B_v).$$
	Choosing $\epsilon$ and $\delta$ to be sufficiently small, we may assume that $B_x\cap B_y=\emptyset$, $B_x \cap C_{uv}=\emptyset$, and $C_{uv}\cap C_{ab}=\emptyset$ for all $x,y\in V(G_1')$ and pairwise non-crossing edges $ab,uv\in E(G_1')$. Let $T^{(1)}$ be a spanning tree of $\tilde{G}_1$ rooted at some $d_1\in \tilde{D}_1$. Using a standard blow-up trick (e.g. see \cite[Lemma~5.5]{HW2021treedensities}), we may replace $T^{(1)}$ by a subdivided star $S^{(1)}$ embedded in the region $(\bigcup_{x\in V(\tilde{G}_1)}B_x)\cup (\bigcup_{uv\in E(\tilde{G}_1)}C_{uv})$ where each root to leaf path in $S^{(1)}$ corresponds to a root to leaf path in $T^{(1)}$; see \cref{fig:SubdividedStar}. 
	\begin{figure}[!htb]
		\begin{center}
			\includegraphics[width=0.7\linewidth]{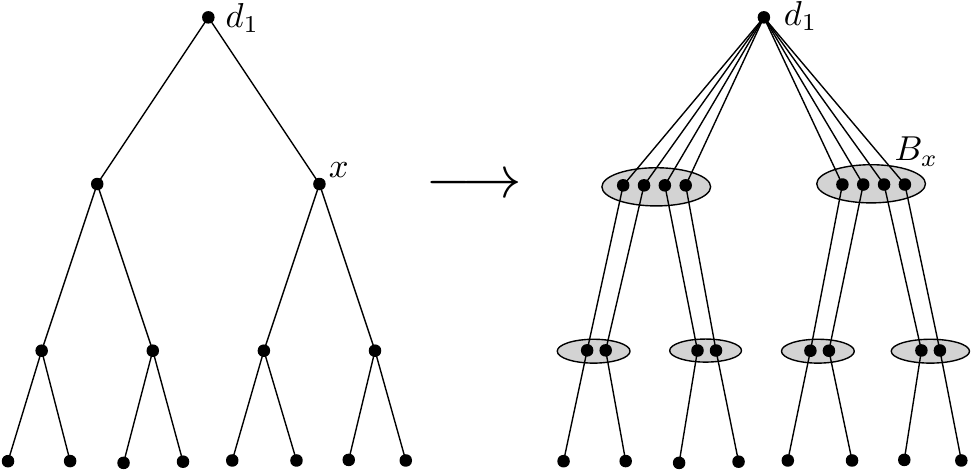}
			\caption{Replacing an embedded tree by a subdivided star.}
			\label{fig:SubdividedStar}
		\end{center}
	\end{figure}
	
	
	Let $G^{(1)}$ be the embedded graph obtained from $G_1'$ by replacing $\tilde{G}_1$ by $S^{(1)}$ then removing the subdivision vertices in $S^{(1)}$; see \cref{fig:FanCrossingEdge}. Charge $u_1$, $v_1$ and $w_1$ to $d_{1}$. By Claim~1, for every crossed edge of the form $d_1z\in E(G^{(1)})$ where $z \in V_1$, there exists an edge $xz\in E_1$ where $x \in \{u_1,v_1,w_1\}$ such that every edge that cross $d_1z$ (in $G^{(1)}$) also crosses $xz$ (in $G^{(0)}$) in the same direction. 
	
	\begin{figure}[!htb]
		\centering
		\includegraphics[width=0.6\linewidth]{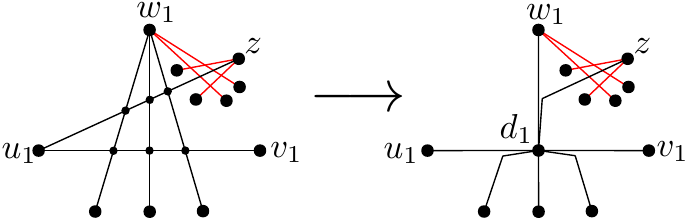}
		\caption{Obtaining $G^{(1)}$ from $G_1'$.}
		\label{fig:FanCrossingEdge}
	\end{figure}
	
	\textbf{Claim 2:} $G^{(1)}$ is simple fan-planar and there is a well-behaved friend assignment of $G^{(1)}$ where $d_1$ is not an assigned friend of any edge.
	\begin{proof}
		 We first show that $G^{(1)}$ is simple fan-planar. Let $e \in E(G^{(1)})$ be a crossed edge and let $A_e$ be the set of edges in $E(G^{(1)})$ that cross $e$. If $\{e\}\cup A_e \sse E(G^{(0)})$, then Configuration I and Configuration II are forbidden and $e$ does not cross any edge at least twice since $G^{(0)}$ is simple fan-planar. So assume that $(\{e\}\cup A_e) \cap (E(G^{(1)})\setminus E(G^{(0)}))\neq \emptyset$. There are two cases to consider.
		
		First, suppose $e=d_{1}z$ for some $z \in V_1$. Since the edges in $E(G^{(1)})\setminus E(G^{(0)})$ do not pairwise cross, it follows that $A_e \sse E(G^{(0)})$ and $e$ does not cross any edge incident to $d_1$. By Claim 1, there is an edge $xz \in E_{1}$ where $x\in \{u_1,v_1,w_1\}$ such that every edge that crosses $e$ (in $G^{(1)}$) also crosses $xz$ (in $G^{(0)}$) in the same direction. Thus the edges that cross $e$ have a common end-vertex and they all cross $e$ in the same direction and $e$ does not cross any edge more than once or any edge incident to $z$.
	
		Second, suppose $e\in E(G^{(1)})\cap E(G^{(0)})$ and $A_e\cap (E(G^{(1)})\setminus E(G^{(0)}))\neq \emptyset$. Let $y$ be the assigned friend of $e$ in $G^{(0)}$. Then $y\in V_1\bs \{u_1,v_1,w_1\}$ since $e$ crosses an edge in $E(G^{(1)})\setminus E(G^{(0)})$ and $E_1$ is maximal. Thus $A_e\cap (E(G^{(1)})\setminus E(G^{(0)}))=\{d_1y\}$ since $y \not \in \{u_1,v_1,w_1\}$. As such, every edge that crosses $e$ (in $G^{(1)}$) is incident to $y$. By Claim 2, there is an edge $xy \in E_{1}$ where $x\in \{u_1,v_1,w_1\}$ such that every edge that crosses $d_1y$ (in $G^{(1)}$) crosses $xy$ (in $G^{(0)}$) in the same direction. Thus Configuration I and II are forbidden and $e$ does not cross any edge more than once since $G^{(0)}$ is simple fan-planar, as required.
		
		We now specify the friend assignment. By \cref{FriendsFan}, $G$ has a well-behaved friend assignment since it is simple fan-planar. Now suppose $d_1$ is the assigned friend of some edge $e$. Then $e$ crosses an edge of the form $d_1y$ where $y \in V_1\bs \{u_1,v_1,w_1\}$. Since every edge that crosses $e$ in $G^{(1)}$ is also incident to $y$, we may modify the friend assignment so that $y$ is the assigned friend of $e$. By doing so for every such edge, we obtain the desired friend assignment.
	\end{proof}
	Observe the following properties of $G^{(1)}$:
	\begin{compactitem}
		\item edges that are uncrossed in $G^{(0)}$ remain uncrossed in $G^{(1)}$;
		\item the edges $u_1d_1,v_1d_1$ and $w_1d_1$ are uncrossed in $G^{(1)}$; and
		\item $G^{(1)}$ has less crossings than $G^{(0)}$.
	\end{compactitem}
	Let $(u_1,v_1,w_1,E_1,d_1,G^{(1)}), (u_2,v_2,w_2,E_2,d_2,G^{(2)}), \dots, (u_m,v_m,w_m,E_m,d_m,G^{(m)})$ be the sequence obtained by iterating the above procedure where $u_iv_i\in E(G^{(i-1)})$ is an arbitrary edge crossed at least twice in $G^{(i-1)}$, $d_i$ is not an assigned friend of any edge in $G^{(j)}$ for all $i\leq j$ and each edge in $G^{(m)}$ is crossed at most once. Then $G^{(m)}$ is $1$-planar. Note that this sequence is well-defined since $G^{(i)}$ has less crossings than $G^{(i-1)}$.
	
	Add a dummy vertex at each crossing point in $G^{(m)}$ and let $H$ be the plane graph obtained. For every edge of the form $d_{i}z\in E(G^{(m)})$ that is crossed by another edge $e$, charge $z$ as well as the end-vertices of $e$ to the dummy vertex at $d_{i}z\cap e$. For crossed edges of the form $e',\tilde{e} \in E(G)\cap E(G^{(m)})$, charge an end-vertex of $e'$ and an end-vertex of $\tilde{e}$ to the dummy vertex at $e'\cap \tilde{e}$. 
	
	Let $D:=V(H)\setminus V(G)$ be the set of dummy vertices in $H$. Observe that every dummy vertex $d \in D$ has at most three vertices charged to it and that these vertices are neighbours of $d$ in $H$. In addition, the original edges in $G$ have the following key properties:

	\textbf{Claim 3:} If $xy\in E(G)$ then either:
	\begin{compactenum}
		\item[Case 1.] $xy\in E(H)$;
		\item[Case 2.] there is a path $P_{xy}=(x,d,y)$ in $H$ where $x$ or $y$ is charged to $d$ and $d\in D$; or
		\item[Case 3.] there is a path $P_{xy}=(x,d_1,d_2,y)$ in $H$ where $x$ is charged to $d_1$ and $y$ is charged to $d_2$ and $d_1,d_2\in D$.
	\end{compactenum}
	\begin{proof}
			First, consider an edge $xy\in E(G^{(m)})\cap E(G)$. If $xy$ is uncrossed in $G^{(m)}$, then $xy \in E(H)$ and we are in Case 1. Otherwise $xy$ is involve in a single crossing in $G^{(m)}$ and so we are in Case 2. So assume $xy\in E_i$ for some $i \in [m]$. Since every edge in $E_i$ is incident to $u_i,v_i$ or $w_i$, we may assume that $x \in \{u_i,v_i,w_i\}$. Thus $x$ is charged to $d_i$ and $d_ix\in E(H)$. If $d_iy\in E(G^{(m)})$ and it is uncrossed, then we are in Case 2. If $d_iy\in E(G^{(m)})$ and it is crossed, then we are in Case 3. So it remains to consider the case when $d_{i}y\in E_{j}$ for some $j\in [i+1,m]$. Since $d_i$ is not the assigned friend of any edge in $G^{(j)}$, it follows that $y\in \{u_j,v_j,w_j\}$. Thus every edge $e$ that crosses $d_iy$ is contained in $E_j$ as $y$ is the assigned friend of $e$. As such, after planarising $E_{j}$ (by replacing it with a star rooted at $d_{j}$), $y$ will be charged to $d_{j}$ and the edges $d_{i}d_{j}$ and $d_jy$ will be uncrossed in $G^{(j)}$, and hence we are in Case 3. 
	\end{proof}
	To finish the proof, we now specify the model $\mu$ for the $1$-shallow minor of $G$ in $H \circ \overline{K_3}$. For each dummy vertex $d \in D$, let $\phi_d: N_d \to \{1,2,3\}$ be an injective function where $N_d$ is the set of vertices that are charged to $d$. For a vertex $u \in V(G)$, let $M_u$ be the set of dummy vertices that $u$ is charged to. For each $u \in V(G)$, let $\mu(u)$ be the subgraph of $H \circ \overline{K_3}$ induced by $\{(u,1)\}\cup \{(d,\phi_d(u)): d \in M_u \}$. 	

	Let $u,v \in V(G)$ be distinct. First $\mu(u)$ is connected and has radius at most $1$ since $ud \in E(H)$ for all $d \in M_u$. Second, if $M_u\cap M_v=\emptyset$, then clearly $\mu(u)$ and $\mu(v)$ are disjoint. Otherwise, if there exists some $d \in M_u \cap M_v$, then $\phi_d(u)\neq \phi_d(v)$ since $\phi_d$ is injective. Thus, $\mu(u)$ and $\mu(v)$ are vertex disjoint. Finally, Claim 3 implies that if $uv \in E(G)$ then $\mu(u)$ and $\mu(v)$ are adjacent. Therefore $\mu$ defines a $1$-shallow model of $G$ in $H \circ \overline{K_3}$. 
\end{proof}

Applying \Cref{FanShallow} with \Cref{GPSTshallow}, we obtain the following product structure theorem for fan-planar graphs.

\begin{thm}\label{FanProductStructure}
	Every fan-planar graph $G$ is contained in $H\boxtimes P\boxtimes K_{81}$ for some graph $H$ with treewidth at most $19$, and thus $G$ has row treewidth at most $1619$.
\end{thm}

In addition, by applying \Cref{FanProductStructure} with \Cref{QueueProduct,NonrepProduct,CenteredColouring}, it follows that for every fan-planar graph $G$,
\begin{compactitem}
	\item $\qn(G)\leq 3\times81\times2^{19}+3\times\floor{81/2}=127402104 $;
	\item $\pi(G)\leq 81\times4^{20}$; and
	\item $\chi_p(G)\leq 81(p+1)\binom{p+19}{19}$.
\end{compactitem}
Note that by using the shallow-minor structure of fan-planar graphs (\cref{FanShallow}) together with \cref{QueueShallow}, Equation~\cref{qnGKl}, and the fact that planar graphs have queue-number at most $42$, it follows that for every fan-planar graph $G$,
$$\qn(G)\leq 2(2\times((2\times 3-1) 42+3-1))^2=359552.$$
Thus we obtain a stronger bound for the queue-number of fan-planar graphs using their shallow-minor structure instead of their product structure.

For colouring numbers, by applying \Cref{FanShallow} in conjunction with \Cref{ColShallowPlanar}, it follows that for every fan-planar graph $G$:
\begin{equation*}
	\scol_s(G)\leq 45s+33 \quad \text{and} \quad	\wcol_s(G) \leq \tbinom{3s+4}{2}(18s+15).
\end{equation*}

For layered treewidth and boxicity, by applying \Cref{FanShallow} with \Cref{ltwShallow,ltwBoxicity}, it follows that $\ltw(G) \leq 45$ and $\boxicity(G)\leq 276$ for every fan-planar graph $G$.

We now show that fan-planar graphs cannot be described by applying a shortcut system to a planar graph. 

\citet[Theorem~3]{BDDMPST2015fan} proved that for every integer $k\geq 1$, there is a fan-planar graph that is not $k$-planar. Essentially the same proof gives the following, stronger result.

\begin{prop}\label{FanNotGap}
	For every integer $k \geq 1$, there is a fan-planar graph that is not $k$-gap planar. 
\end{prop}

\begin{wrapfigure}{R}{0.48\textwidth}
		\centering
		\includegraphics[angle=90,width=0.8\linewidth,height=2.8cm]{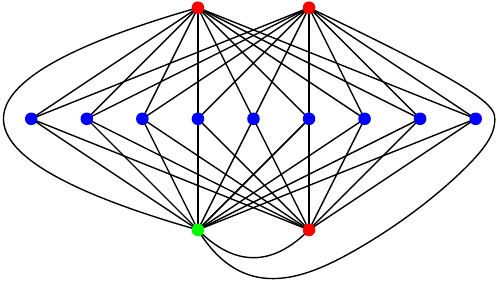}
		\caption{A fan-planar embedding\\ of $K_{1,3,h}$.}
		\label{fig:FanPlanarExample}
		\vspace*{-0.7cm}
	\end{wrapfigure}

\begin{proof}
	For each $h \in \N$, the complete tripartite graph $K_{1,3,h}$ is fan-planar and $|E((K_{1,3,h}))|=4h+3$; see \Cref{fig:FanPlanarExample}. Moreover, it is known \cite[Theorem~1]{asano1986crossing} that $\cross(K_{1,3,h})=\Omega(h^2)=\Omega(|E((K_{1,3,h}))|^2)$. Thus, this graph family has super-linear crossing number so it is not $k$-gap planar. 
\end{proof}

\Cref{ShortcutGapPlanar,FanNotGap} imply that fan-planar graphs cannot be described by shortcut systems applied to a planar graph. This highlights the value and power of shallow minors.

As a converse to \cref{FanNotGap}, \citet[Theorem~4]{BDDMPST2015fan} showed that there exists a $2$-planar graph that is not fan-planar. Since every $2$-planar graph is $1$-gap planar graph \cite[Theorem~7]{BBCEEGHKMRT2018gap}, we have the following.
\begin{prop}\label{GapNotFan}
	For every integer $k\geq 1$, there is a $k$-gap planar graph that is not fan-planar.
\end{prop}

\Cref{FanNotGap,GapNotFan} demonstrates that for any integer $k\geq 1$, the class of fan-planar graphs and $k$-gap planar graphs are incomparable.\footnote{Graph classes $\G_1$ and $\G_2$ are \defn{incomparable} if $\G_1 \bs \G_2\neq \emptyset$ and $\G_2 \bs \G_1 \neq \emptyset$.} This resolves an open problem of \citet[Problem~12]{DLM2019beyond} who asked for the relationship between $k$-gap planar graphs and fan-planar graphs.

\subsection{\normalsize\boldmath $k$-Fan-Bundle Planar Graph}\label{SectionkFanBundle}

In a $k$-fan-bundle embedding of a graph in the plane the edges of a fan may be bundled together at their end-vertices and crossings between bundles are allowed as long as each bundle is crossed by at most $k$ other bundles. More formally, in a $k$-fan-bundle planar embedding of a graph $G$, each edge has three parts; the first and the last parts are \defn{fan-bundles}, which may be shared by several edges in a fan, while the middle part is unbundled. Each fan-bundle can cross at most $k$ other fan-bundles, while the unbundled parts are crossing-free. A vertex $u\in V(G)$ can be incident to more than one bundle. Let $B_u$ be one such bundle. We say that $B_u$ is \defn{anchored} at $u$ which is the \defn{origin} of $B_u$. The endpoint of $B_u$ that is different from $u$ is the \defn{terminal} of $B_u$. A graph is \defn{$k$-fan-bundle planar} if it admits a \defn{$k$-fan-bundle-planar embedding}. $k$-fan bundle planar graphs were introduced by \citet{ABKKS2018fanbundle}.

\begin{lem}\label{kFanBundleShallow}
	Every $k$-fan-bundle planar graph $G$ is a $(k+1)$-shallow minor of $H \circ \overline{K_2}$ for some planar graph $H$.
\end{lem} 

\begin{wrapfigure}{R}{0.33\textwidth}
		\hfill\includegraphics[width=0.75\linewidth]{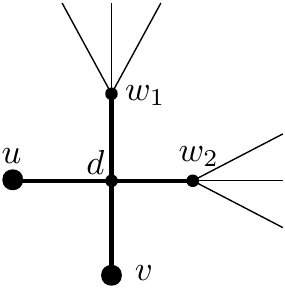}
		\caption{Planarising a\\ $1$-fan-bundle planar graph.}
		\label{fig:fanbundle}
	\end{wrapfigure}

\beginproof \ 
	Begin with a $k$-fan-bundle-planar embedding of $G$. For each bundle $B_u^{(i)}$, add a dummy vertex at the terminal of $B_u^{(i)}$ to obtain a $k$-planar graph $H'$. Let $W$ be the set of dummy vertices added. Replace each crossing in $H'$ by a dummy vertex to obtain a plane graph $H$ (see \Cref{fig:fanbundle}) and let $D$ be the new dummy vertices that are added at this step. In doing so, each bundle $B_u^{(i)}$ is replaced by a $(u,w)$-path $P^{(i)}_u$ on at most $k+2$ vertices for some dummy vertex $w\in W$. For each $d \in D$, let $\phi_d: \{u,v\} \to \{1,2\}$ be an injective function where $u$ and $v$ are the origins of the bundles that cross at $d$. For each $u \in V(G)$, let $W_u:=(\bigcup_{i\in [c_u]}V(P^{(i)}_u))\cap W$ and $D_u:=(\bigcup_{i\in [c_u]}V(P^{(i)}_u))\cap D$ where $c_u$ is the number of bundles anchored at $u$. For each $u \in V(G)$, let $\mu(u)$ be the subgraph of $H \circ \overline{K_2}$ induced by $\{(u,1)\}\cup \{(w,1): w \in W_u \}\cup \{(d,\phi_w(d)): d \in D_u \}$. 
	
	We claim that $\mu$ is a $(k+1)$-shallow model of $G$ in $H \circ \overline{K_2}.$ 
	Let $u,v \in V(G)$ be distinct. First, $\mu(u)$ is connected with radius at most $k+1$ as it is the union of paths on at most $k+2$ vertices that share a common end-vertex. Now if no bundle anchored at $u$ crosses a bundle anchored at $v$, then $\mu(u)$ and $\mu(v)$ are clearly disjoint. Otherwise, if there is a bundle $B_{u}^{(i)}$ that crosses a bundle $B_{v}^{(i')}$, then there is some $d \in D_u \cap D_v$. In which case, $\phi_d(u)\neq \phi_d(v)$ since $\phi_d$ is injective and so $\mu(u)$ and $\mu(v)$ are vertex disjoint. Now if $uv \in E(G)$, then there is a bundle $B_u^{(i)}$ anchored at $u$ and a bundle $B_v^{(i')}$ anchored at $v$ that are adjacent. Let $w_1,w_2\in W$ be the dummy vertices that are respectively added to the terminals of $B_u^{(i)}$ and $B_v^{(i')}$. Then $(w_1,1)\in \mu(u)$ and $(w_2,1)\in \mu(v)$ so $(w_1,1)(w_2,1)\in E(H \circ \overline{K_2})$ and thus $\mu(u)$ and $\mu(v)$ are adjacent, as required.	\qed
	
Applying \Cref{kFanBundleShallow} with \Cref{GPSTshallow}, we obtain the following product structure theorem for $k$-fan-bundle planar graphs.

\begin{thm}\label{fanbundleGPST}
	Every $k$-fan-bundle planar graph $G$ is contained in $H\boxtimes P\boxtimes K_{6(2k+3)^2}$ for some graph $H$ with treewidth at most $\binom{2k+6}{3}-1$, and thus $G$ has row treewidth at most $\binom{2k+6}{3}6(2k+3)^2-1$.
\end{thm}

Applying \Cref{fanbundleGPST} with \Cref{QueueProduct,NonrepProduct,CenteredColouring}, it follows that every $k$-fan-bundle planar graph $G$ has:
\begin{compactitem}
	\item $\qn(G)\leq 18(2k+3)^22^{\binom{2k+6}{3}-1}+9(2k+3)^2$;
	\item $\pi(G)\leq 6(2k+3)^24^{\binom{2k+6}{3}}$; and
	\item $\chi_p(G)\leq 6(2k+3)^2(p+1)\binom{p+\binom{2k+6}{3}-1}{\binom{2k+6}{3}-1}$.
\end{compactitem}

Note that by using the shallow-minor structure of $k$-fan-bundle planar graph (\cref{kFanBundleShallow}) together with \cref{QueueShallow,qnGKl}, and the fact that planar graphs have queue-number at most $42$, it follows that for every $k$-fan-bundle planar graph $G$, 
$$\qn(G)\leq 2(k+1)(2\times((2\times 2-1) 42+2-1))^{2(k+1)}= 2(k+1)(170)^{2(k+1)}.$$
Thus we obtain a stronger bound for the queue-number of $k$-fan-bundle planar graph using their shallow-minor structure instead of their product structure.

For colouring numbers, \Cref{kFanBundleShallow} and \Cref{ColShallowPlanar} imply that for every $k$-fan-bundle planar graph $G$:
\begin{align*}
	\scol_s(G) &\leq (20k+30)s+20k+22\quad\text{and}\\
	\wcol_s(G) &\leq \tbinom{(2k+3)s+2k+4}{2}((8(k+3)s+8k+10).
\end{align*}

For layered treewidth and boxicity, it follows from \Cref{kFanBundleShallow,ltwShallow,ltwBoxicity}, that $\ltw(G) \leq 24k+25$ and $\boxicity(G)\leq 144k+154$ for every $k$-fan-bundle planar graph $G$.

We raise the following open problems concerning $k$-fan-bundle graphs. \citet{ABKKS2018fanbundle} showed that $1$-fan-bundle graphs are incomparable with $2$-planar graphs. What is the relationship between $1$-fan-bundle graphs and $k$-planar graphs? Does there exist an integer $k\geq 1$ such that every $1$-fan bundle graph is $k$-planar? Or more weaker, does there exists an integer $k\geq 1$ such that every $1$-fan bundle graph is $k$-gap planar? 

\section{\large Lower Bounds}\label{SectionLowerBounds}
Having described several beyond planar graph classes as shallow minors of the strong product of a planar graph with a small complete graph, we now give examples of classes that cannot be described in this manner.

\subsection{\normalsize\boldmath $k$-Gap Planar}\label{SectionGap}

Recall that a graph $G$ is $k$-gap planar if it is isomorphic to an embedded graph where each crossing is charged to one of the two edges involved and each edge has at most $k$ crossings charged to it. This class of graphs has been implicitly studied for some time \cite{OOW2019defective,EG2017road}. The language of $k$-gap planar graphs was introduced by \citet{BBCEEGHKMRT2018gap}. We show that $k$-gap planar graphs have unbounded row treewidth and thus cannot be described as a shallow minor of the strong product of a planar graph with a small complete graph, even with $k=1$. This result is of particular interest since $k$-gap planar graphs have polynomial expansion \cite[Theorem~6.9]{EG2017road}. We in fact show a stronger result in terms of local treewidth.

\citet{eppstein2000diameter} introduced the following definition under the guise of the `treewidth-diameter' property. A graph class $\G$ has \defn{bounded local treewidth} if there is a function $f$ such that for every graph $G \in \G$, for every vertex $v\in V(G)$ and for every integer $r\geq 0$, the subgraph of $G$ induced by the vertices at distance at most $r$ from $v$ has treewidth at most $f(r)$. If $f(r)$ is linear, then $\G$ has \defn{linear local treewidth}. For example, every class with bounded layered treewidth (which includes graph classes with bounded row treewidth) has linear local treewidth (\cite[Lemma~6]{dujmovic2017layered}). We show that $1$-gap planar graphs do not have linear local treewidth which implies that they have unbounded row treewidth. In fact, we prove an exponential lower bound on the local treewidth.

\begin{thm}\label{GapPlanar}
	There exists a constant $c>0$ such that for infinitely many integers $r\geq 1$, there is a $1$-gap planar graph $G$ with radius at most $r$ and $\tw(G)\geq 2^{cr}$.
\end{thm}
\Cref{GapPlanar} is implied by the following lemma by setting $n=2$. 

\begin{lem}\label{GapPlanarLemmaNew}
	For all integers $n,k\geq 1$, there exists a $1$-gap planar graph with treewidth at least $n^{k+1}+1$ and radius at most $(2k+1)n+\ceil{\frac{k}{2}}+1$.
\end{lem}
\begin{proof}
	Let $G_{\infty}$ be the following infinite graph with $V(G_{\infty}):=\mathbb{Z}^2$. For ${\ell} \in [0,k]$, let $A^{(\ell)}:=\{(an^{\ell},bn^{\ell}):a,b \in \mathbb{Z}\}$. Add the \defn{level-${\ell}$ horizontal edges} $(an^{\ell},bn^{\ell})((a+1)n^{\ell},bn^{\ell})$ and the \defn{level-${\ell}$ vertical edges} $(an^{\ell},bn^{\ell})(an^{\ell},(b+1)n^{\ell})$ to $G_{\infty}$ for all $a,b\in \mathbb{Z}$. 
	
	\begin{figure}[h]
		\centerline{\includegraphics[width=0.83\textwidth]{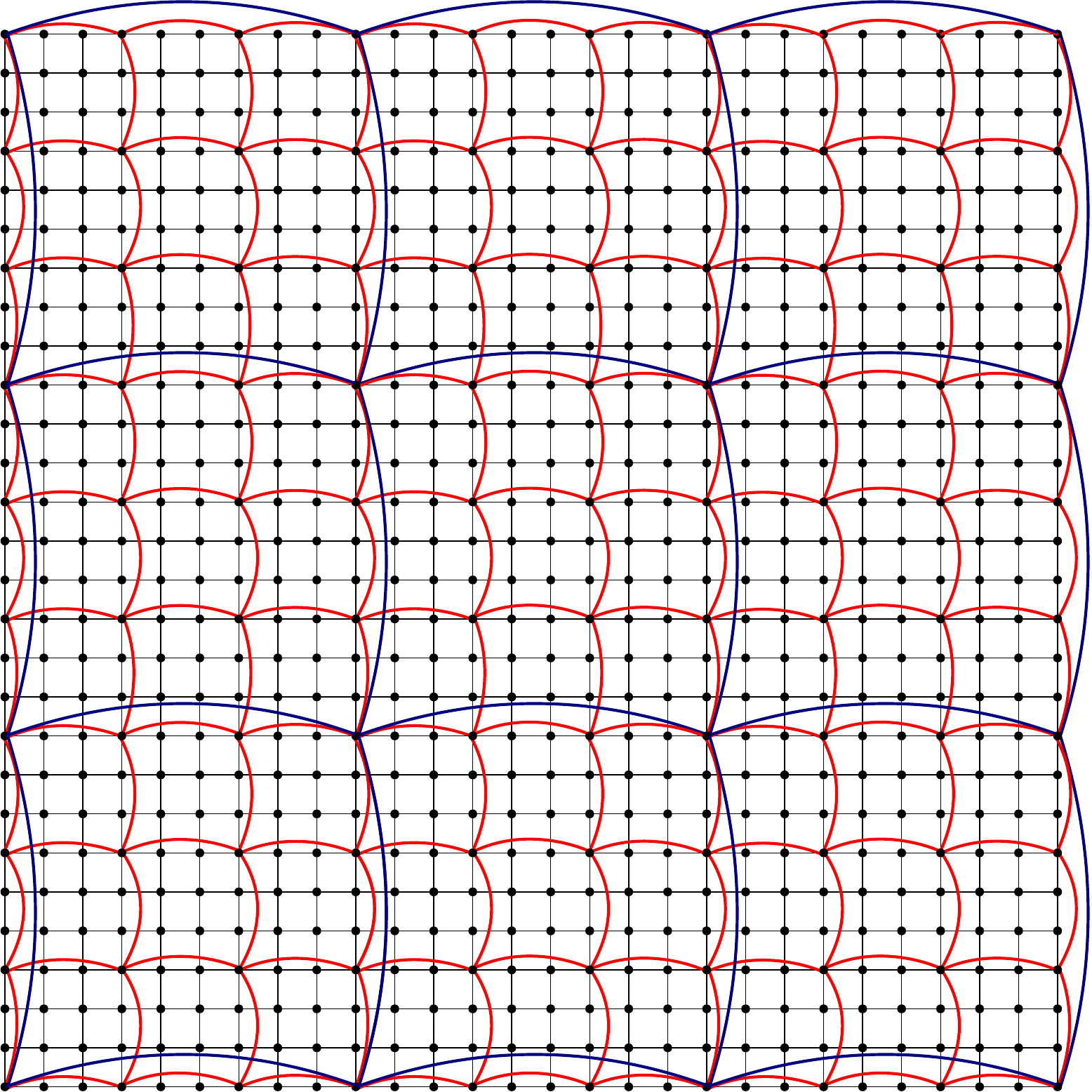}}
		\caption{Construction in the proof of \Cref{GapPlanarLemmaNew} with $k=2$ and $n=3$.}
		\label{fig:GapPlanar}
	\end{figure}
	
	We claim that $G_{\infty}$ is $k$-gap planar (see \Cref{fig:GapPlanar}). Embed each vertex $(i,j)\in V(G_{\infty})$ at $(i,j)\in \R^2$. Draw the level-$0$ edges as straight line segments between their end-vertices. For ${\ell} =1,2\dots,k$ and $a,b\in \mathbb{Z}$, draw the edge $(an^{\ell},bn^{\ell})((a+1)n^{\ell},bn^{\ell})$ as a curve between the $y=bn^{\ell}$ line and $y=bn^{\ell}+1$ line such that it does not cross any horizontal edges or level-$\ell$ edges. Similarly, draw the edge $(an^{\ell},bn^{\ell})(an^{\ell},(b+1)n^{\ell})$ as a curve between the $x=an^{\ell}$ line and $x=an^{\ell}+1$ line such that it does not cross any vertical edges or level-$\ell$ edge. This can be done so that for every $j,\ell\in [0,k]$ where $j <\ell$, every level-$j$ horizontal edge is crossed by at most one vertical level-${\ell}$ edge and every level-$j$ vertical edge is crossed by at most one horizontal level-${\ell}$ edge. Charge the crossing to the edge with the lower level. Thus every edge in $G_{\infty}$ has at most $k$ crossings charged to it and so $G_{\infty}$ is $k$-gap planar.
	
	Let $G_n:=G_{\infty}[[0,n^{k+1}]\times [0,n^{k+1}]]$. Then $G_n$ is $k$-gap planar since it is a subgraph of $G_{\infty}$. Moreover, $G_n$ has treewidth at least $n^{k+1}+1$ since it contains the $(n^{k+1}+1) \times (n^{k+1}+1)$-grid as a subgraph \cite{bodlaender1998partial}. 
	We now subdivide edges in $G_n$ to obtain a $1$-gap planar graph $\tilde{G}_n$. Replace each edge $uv \in E(G_n)$ that has $c>1$ crossings charged to it by a $(u,v)$-path, $P_{uv}$, on $c+1$ vertices such that the internal vertices of $P_{uv}$ are between consecutive crossing points that are charged to $uv$. Let $G_n$ be the $1$-gap planar graph obtained from this subdivision procedure. Observe that for every $\ell \in [0,k-1]$, every horizontal subdivided edge between two vertices in $A^{(\ell)}\cap V(\tilde{G}_n)$ is of the form $(an^{\ell+1},bn^{\ell})((an+1)n^{\ell},bn^{\ell})$ and every vertical subdivided edge is of the form $(an^{\ell},bn^{\ell+1})(an^{\ell},(bn+1)n^{\ell})$. Moreover, the level-$k$ horizontal and vertical edges are not subdivided as no crossings are charged to them.
	
	Now $\tilde{G}_n$ has treewidth at least $n^{k+1}+1$ since it is a subdivision of $G_n$. So it remains to show that $\tilde{G}_n$ has radius at most $(2k+1)n+\ceil{\frac{k}{2}}+1$. First, every internal vertex of a subdivided path $P_{uv}$ has distance at most $\ceil{\frac{k}{2}}$ to $u$ or $v$ where $u,v\in V(\tilde{G}_n)\cap A^{(0)}$. Let $\ell \in [k]$. Observe that every vertex in $A^{(\ell-1)}\cap V(\tilde{G}_n)$ is of the form $((an+i)n^{\ell-1},(bn+j)n^{\ell-1})$ where $a,b\geq 0$ are integers and $i,j \in [n]$. Thus $\tilde{G}_n$ contains the following path from $((an+i)n^{\ell-1},(bn+j)n^{\ell-1})\in A^{(\ell-1)}\cap V(\tilde{G}_n)$ to $((a+1)n^{\ell},(b+1)n^{\ell}) \in A^{(\ell)}\cap V(\tilde{G}_n)$ which avoids subdivided paths and has length at most $2n$:
	\begin{align*}
		&((an+i)n^{\ell-1},(bn+j)n^{\ell-1}),((an+i+1)n^{\ell-1},(bn+j)n^{\ell-1}),\dots, ((a+1)n^{\ell},(bn+j)n^{\ell-1}),\\
		&((a+1)n^{\ell},(bn+j+1)n^{\ell-1}),\dots, ((a+1)n^{\ell},(b+1)n^{\ell}).
	\end{align*}
	As such, for each $v\in V(\tilde{G}_n)$, there is a path of length at most $2kn+\ceil{\frac{k}{2}}$ from $v$ to some vertex in $A^{(k)}\cap V(\tilde{G}_n)$. Finally, the induced subgraph $\tilde{G}_n[A^{(k)}]$ has radius at most $n+1$ since it is isomorphic to the $(n+1) \times (n+1)$ grid. Therefore $\tilde{G}_n$ has radius at most $(2k+1)n+\ceil{\frac{k}{2}}+1$.
\end{proof}

We conclude this subsection with the following open problem: do $k$-gap planar graphs have bounded local treewidth?

\subsection{\normalsize\boldmath RAC, Fan-Crossing Free, and $k$-Quasi-Planar Graphs}\label{SectionSomewhereDense}

Recall that an {embedded graph $G$} is:
\begin{compactitem}
	\item {$k$-quasi-planar }if every set of $k$ edges do not mutually cross;
	\item {fan-crossing free} if for each edge $e \in E(G)$, the edges that cross $e$ form a matching; or 
	\item {right angle crossing }(RAC) if each edge is drawn as a straight line segment and all crossings are at right angles.
\end{compactitem}

A graph is {$k$-quasi planar}, {fan-crossing free}, or {RAC} if the graph is respectively isomorphic to an embedded graph that is $k$-quasi planar, fan-crossing free, or RAC. 
Recall that a graph class $\G$ is {somewhere dense} if there exists an integer $r\geq 0$ such that every graph $H$ is an $r$-shallow minor of some graph $G \in \G$ (see \cite{nevsetvril2012sparsity}). Drawing on known results, \citet{brandenburg2021fancrossingfree} proved that for every graph $G$, 
the $3$-subdivision of $G$ is a RAC graph;
the $2$-subdivision of $G$ is fan-crossing free; and
the $1$-subdivision of $G$ is $3$-quasi-planar.
This implies that each of these classes are somewhere dense. Note that \citet{eppstein2015shallow} previously observed that the class of RAC graphs is somewhere dense. Since graph classes with bounded row treewidth have polynomial expansion \cite{dvovrak2020notes}, these beyond planar classes have unbounded row treewidth. Thus, they cannot be described as a shallow minor of the strong product of a planar graph with a small complete graph. Additionally, \citet{brandenburg2021fancrossingfree} asked what is the queue-number and stack-number of these classes. Since graph classes with bounded stack-number or bounded queue-number have bounded expansion \cite{NOW2012examples}, we conclude that RAC, fan-crossing free, and $k$-quasi planar graphs (where $k\geq 3$) have unbounded stack-number and unbounded queue-number. 	

\fontsize{9.5}{10.5} 
\selectfont 
\let\oldthebibliography=\thebibliography
\let\endoldthebibliography=\endthebibliography
\renewenvironment{thebibliography}[1]{%
	\begin{oldthebibliography}{#1}%
		\setlength{\parskip}{0.2ex}%
		\setlength{\itemsep}{0.2ex}%
	}{\end{oldthebibliography}}
\bibliographystyle{DavidNatbibStyle}
\bibliography{main.bbl}	

\begin{thebibliography}{83}
\providecommand{\natexlab}[1]{#1}
\providecommand{\msn}[1]{MR:\,\href{http://www.ams.org/mathscinet-getitem?mr=MR{#1}}{#1}}
\providecommand{\ZBL}[1]{Zbl:\,\href{https://www.zentralblatt-math.org/zmath/en/search/?q=an:#1}{#1}}
\providecommand{\url}[1]{\texttt{#1}}
\providecommand{\urlprefix}{}
\expandafter\ifx\csname urlstyle\endcsname\relax
  \providecommand{\doi}[1]{doi:\discretionary{}{}{}#1}\else
  \providecommand{\doi}{doi:\discretionary{}{}{}\begingroup
  \urlstyle{rm}\Url}\fi

\bibitem[{Ajtai et~al.(1982)Ajtai, Chv{\'a}tal, Newborn, and
  Szemer{\'e}di}]{ACNS1982crossing}
\textsc{Mikl{\'o}s Ajtai, Va{\v{s}}ek Chv{\'a}tal, Monroe~M Newborn, and Endre
  Szemer{\'e}di}.
\newblock Crossing-free subgraphs.
\newblock In \emph{Theory and practice of combinatorics}, vol.~60 of
  \emph{North-Holland Math. Stud.}, pp. 9--12. North-Holland, 1982.

\bibitem[{Alon et~al.(2002)Alon, Grytczuk, Ha{\l}uszczak, and
  Riordan}]{AGHR2002nonrepetitive}
\textsc{Noga Alon, Jaros{\l}aw Grytczuk, Mariusz Ha{\l}uszczak, and Oliver
  Riordan}.
\newblock \href{https://doi.org/10.1002/rsa.10057}{Nonrepetitive colorings of
  graphs}.
\newblock \emph{Random Structures Algorithms}, 21(3-4):336--346, 2002.

\bibitem[{Angelini et~al.(2018)Angelini, Bekos, Kaufmann, Kindermann, and
  Schneck}]{ABKKS2018fanbundle}
\textsc{Patrizio Angelini, Michael~A. Bekos, Michael Kaufmann, Philipp
  Kindermann, and Thomas Schneck}.
\newblock \href{https://doi.org/10.1016/j.tcs.2018.03.005}{1-fan-bundle-planar
  drawings of graphs}.
\newblock \emph{Theoret. Comput. Sci.}, 723:23--50, 2018.

\bibitem[{Asano(1986)}]{asano1986crossing}
\textsc{Kouhei Asano}.
\newblock \href{https://doi.org/10.1002/jgt.3190100102}{The crossing number of
  {$K\sb{1,3,n}$} and {$K\sb{2,3,n}$}}.
\newblock \emph{J. Graph Theory}, 10(1):1--8, 1986.

\bibitem[{Bae et~al.(2018)Bae, Baffier, Chun, Eades, Eickmeyer, Grilli, Hong,
  Korman, Montecchiani, Rutter, and T\'{o}th}]{BBCEEGHKMRT2018gap}
\textsc{Sang~Won Bae, Jean-Francois Baffier, Jinhee Chun, Peter Eades, Kord
  Eickmeyer, Luca Grilli, Seok-Hee Hong, Matias Korman, Fabrizio Montecchiani,
  Ignaz Rutter, and Csaba~D. T\'{o}th}.
\newblock \href{https://doi.org/10.1016/j.tcs.2018.05.029}{Gap-planar graphs}.
\newblock \emph{Theoret. Comput. Sci.}, 745:36--52, 2018.

\bibitem[{Bekos et~al.(2021)Bekos, Gronemann, and Raftopoulou}]{BGR2021queue}
\textsc{Michael~A Bekos, Martin Gronemann, and Chrysanthi~N Raftopoulou}.
\newblock \href{http://arxiv.org/abs/2106.08003}{On the queue number of planar
  graphs}.
\newblock In \emph{Proc. of the 29th {I}nternational {S}ymposium on {G}raph
  {D}rawing and {N}etwork {V}isualization ({GD} 2021)}. 2021.
\newblock arXiv:2106.08003.

\bibitem[{Binucci et~al.(2015)Binucci, Di~Giacomo, Didimo, Montecchiani,
  Patrignani, Symvonis, and Tollis}]{BDDMPST2015fan}
\textsc{Carla Binucci, Emilio Di~Giacomo, Walter Didimo, Fabrizio Montecchiani,
  Maurizio Patrignani, Antonios Symvonis, and Ioannis~G. Tollis}.
\newblock \href{https://doi.org/10.1016/j.tcs.2015.04.020}{Fan-planarity:
  properties and complexity}.
\newblock \emph{Theoret. Comput. Sci.}, 589:76--86, 2015.

\bibitem[{Bodlaender(1993)}]{bodlaender1993tourist}
\textsc{Hans~L. Bodlaender}.
\newblock
  \href{https://www.proquest.com/scholarly-journals/tourist-guide-through-treewidth/docview/2384586349/se-2?accountid=12528}{A
  tourist guide through treewidth}.
\newblock \emph{Acta Cybernet.}, 11(1-2):1--21, 1993.

\bibitem[{Bodlaender(1998)}]{bodlaender1998partial}
\textsc{Hans~L. Bodlaender}.
\newblock \href{https://doi.org/10.1016/S0304-3975(97)00228-4}{A partial
  {$k$}-arboretum of graphs with bounded treewidth}.
\newblock \emph{Theoret. Comput. Sci.}, 209(1-2):1--45, 1998.

\bibitem[{Bodlaender et~al.(1995)Bodlaender, Gilbert, Hafsteinsson, and
  Kloks}]{BGHK1995approximating}
\textsc{Hans~L. Bodlaender, John~R. Gilbert, Hj\'{a}lmt\'{y}r Hafsteinsson, and
  Ton Kloks}.
\newblock \href{https://doi.org/10.1006/jagm.1995.1009}{Approximating
  treewidth, pathwidth, frontsize, and shortest elimination tree}.
\newblock \emph{J. Algorithms}, 18(2):238--255, 1995.

\bibitem[{Bonamy et~al.(2020)Bonamy, Gavoille, and
  Pilipczuk}]{bonamy2020shorter}
\textsc{Marthe Bonamy, Cyril Gavoille, and Micha{\l} Pilipczuk}.
\newblock \href{https://doi.org/10.1137/1.9781611975994.27}{Shorter labeling
  schemes for planar graphs}.
\newblock In \emph{Proc. of 14th Annual ACM-SIAM Symposium on Discrete
  Algorithms}, pp. 446--462. SIAM, 2020.

\bibitem[{Bonnet et~al.(2021)Bonnet, Geniet, Kim, Thomass\'{e}, and
  Watrigant}]{BGKTW2021twinsmall}
\textsc{\'{E}douard Bonnet, Colin Geniet, Eun~Jung Kim, St\'{e}phan
  Thomass\'{e}, and R\'{e}mi Watrigant}.
\newblock \href{https://doi.org/10.1137/1.9781611976465.118}{Twin-width {II}:
  small classes}.
\newblock In \emph{Proc. of the 2021 ACM-SIAM {S}ymposium on {D}iscrete
  {A}lgorithms (SODA)}, pp. 1977--1996. 2021.

\bibitem[{Bonnet et~al.(2022)Bonnet, Kwon, and Wood}]{BKW2022bandwidth}
\textsc{{\'E}douard Bonnet, O-joung Kwon, and David~R. Wood}.
\newblock \href{http://arxiv.org/abs/2202.11858}{Reduced bandwidth: a
  qualitative strengthening of twin-width in minor-closed classes (and
  beyond)}.
\newblock 2022.
\newblock arXiv:2202.11858.

\bibitem[{Bose et~al.(2022)Bose, Dujmovi\'c, Javarsineh, Morin, and
  Wood}]{BDJMW2021rowtreewidth}
\textsc{Prosenjit Bose, Vida Dujmovi\'c, Mehrnoosh Javarsineh, Pat Morin, and
  David~R. Wood}.
\newblock \href{https://doi.org/10.46298/dmtcs.7458}{Separating layered
  treewidth and row treewidth}.
\newblock \emph{Discret. Math. Theor. Comput. Sci.}, 24(1):\#18, 2022.

\bibitem[{Brandenburg(2021)}]{brandenburg2021fancrossingfree}
\textsc{Franz~J. Brandenburg}.
\newblock \href{https://doi.org/10.1016/j.tcs.2021.03.031}{Fan-crossing free
  graphs and their relationship to other beyond-planar graphs}.
\newblock \emph{Theoret. Comput. Sci.}, 867:85--100, 2021.

\bibitem[{Chen and Schelp(1993)}]{CS1993ramsey}
\textsc{Guantao Chen and R.~H. Schelp}.
\newblock \href{https://doi.org/10.1006/jctb.1993.1012}{Graphs with linearly
  bounded {R}amsey numbers}.
\newblock \emph{J. Combin. Theory Ser. B}, 57(1):138--149, 1993.

\bibitem[{Cheong et~al.(2015)Cheong, Har-Peled, Kim, and
  Kim}]{CHKK2015fancrossing}
\textsc{Otfried Cheong, Sariel Har-Peled, Heuna Kim, and Hyo-Sil Kim}.
\newblock \href{https://doi.org/10.1007/s00453-014-9935-z}{On the number of
  edges of fan-crossing free graphs}.
\newblock \emph{Algorithmica}, 73(4):673--695, 2015.

\bibitem[{Debski et~al.(2021)Debski, Felsner, Micek, and
  Schr{\"o}der}]{debski2020improved}
\textsc{Micha{\l} Debski, Stefan Felsner, Piotr Micek, and Felix Schr{\"o}der}.
\newblock \href{https://doi.org/10.19086/aic.27351}{Improved bounds for
  centered colorings}.
\newblock \emph{Adv. Comb.}, (8), 2021.

\bibitem[{Di~Giacomo et~al.(2019)Di~Giacomo, Lenhart, Liotta, Randolph, and
  Tappini}]{DLLRT2019hybrid}
\textsc{Emilio Di~Giacomo, William~J. Lenhart, Giuseppe Liotta, Timothy~W.
  Randolph, and Alessandra Tappini}.
\newblock
  \href{https://doi.org/10.1007/978-3-030-10564-8_12}{{$(k,p)$}-planarity: a
  relaxation of hybrid planarity}.
\newblock In \emph{W{ALCOM}: algorithms and computation}, vol. 11355 of
  \emph{Lecture Notes in Comput. Sci.}, pp. 148--159. Springer, 2019.

\bibitem[{Didimo et~al.(2011)Didimo, Eades, and Liotta}]{DEL2011RAC}
\textsc{Walter Didimo, Peter Eades, and Giuseppe Liotta}.
\newblock \href{https://doi.org/10.1016/j.tcs.2011.05.025}{Drawing graphs with
  right angle crossings}.
\newblock \emph{Theoret. Comput. Sci.}, 412(39):5156--5166, 2011.

\bibitem[{Didimo et~al.(2019)Didimo, Liotta, and Montecchiani}]{DLM2019beyond}
\textsc{Walter Didimo, Giuseppe Liotta, and Fabrizio Montecchiani}.
\newblock \href{https://doi.org/10.1145/3301281}{A survey on graph drawing
  beyond planarity}.
\newblock \emph{ACM Comput. Surv.}, 52(1), 2019.

\bibitem[{Diestel(2017)}]{diestel2017graphtheory}
\textsc{Reinhard Diestel}.
\newblock \href{https://doi.org/10.1007/978-3-662-53622-3}{Graph theory}, vol.
  173 of \emph{Graduate Texts in Mathematics}.
\newblock Springer, 5th edn., 2017.

\bibitem[{Distel et~al.(2021)Distel, Hickingbotham, Huynh, and Wood}]{DHHW21}
\textsc{Marc Distel, Robert Hickingbotham, Tony Huynh, and David~R. Wood}.
\newblock \href{http://arxiv.org/abs/2112.10025}{Improved product structure for
  graphs on surfaces}.
\newblock 2021.
\newblock arXiv:2112.10025.

\bibitem[{Dujmović et~al.(2020)Dujmović, Morin, and Wood}]{DMW20}
\textsc{Vida Dujmović, Pat Morin, and David~R. Wood}.
\newblock \href{http://arxiv.org/abs/1907.05168}{Graph product structure for
  non-minor-closed classes}.
\newblock 2020.
\newblock arXiv:1907.05168.

\bibitem[{Dujmovi{\'c} et~al.(2021)Dujmovi{\'c}, Eppstein, Hickingbotham,
  Morin, and Wood}]{dujmovic2020stack}
\textsc{Vida Dujmovi{\'c}, David Eppstein, Robert Hickingbotham, Pat Morin, and
  David~R Wood}.
\newblock
  \href{https://doi.org/https://doi.org/10.1007/s00493-021-4585-7}{Stack-number
  is not bounded by queue-number}.
\newblock \emph{Combinatorica}, 2021.

\bibitem[{Dujmovi\'{c} et~al.(2021)Dujmovi\'{c}, Esperet, Gavoille, Joret,
  Micek, and Morin}]{DEGJMM2020adjacency}
\textsc{Vida Dujmovi\'{c}, Louis Esperet, Cyril Gavoille, Gwena\"{e}l Joret,
  Piotr Micek, and Pat Morin}.
\newblock \href{https://doi.org/10.1145/3477542}{Adjacency labelling for planar
  graphs (and beyond)}.
\newblock \emph{J. ACM}, 68(6), 2021.

\bibitem[{Dujmovi\'{c} et~al.(2020{\natexlab{a}})Dujmovi\'{c}, Esperet, Joret,
  Walczak, and Wood}]{DEJWW2020nonrepetitive}
\textsc{Vida Dujmovi\'{c}, Louis Esperet, Gwena\"{e}l Joret, Bartosz Walczak,
  and David~R. Wood}.
\newblock \href{https://doi.org/10.19086/aic.12100}{Planar graphs have bounded
  nonrepetitive chromatic number}.
\newblock \emph{Adv. Comb.}, (5), 2020{\natexlab{a}}.

\bibitem[{Dujmovi\'{c} et~al.(2020{\natexlab{b}})Dujmovi\'{c}, Joret, Micek,
  Morin, Ueckerdt, and Wood}]{DJMMUW20}
\textsc{Vida Dujmovi\'{c}, Gwena\"{e}l Joret, Piotr Micek, Pat Morin, Torsten
  Ueckerdt, and David~R. Wood}.
\newblock \href{https://doi.org/10.1145/3385731}{Planar graphs have bounded
  queue-number}.
\newblock \emph{J. ACM}, 67(4):Art. 22, 2020{\natexlab{b}}.

\bibitem[{Dujmovi\'{c} et~al.(2005)Dujmovi\'{c}, Morin, and
  Wood}]{DMW2005layouts}
\textsc{Vida Dujmovi\'{c}, Pat Morin, and David~R. Wood}.
\newblock \href{https://doi.org/10.1137/S0097539702416141}{Layout of graphs
  with bounded tree-width}.
\newblock \emph{SIAM J. Comput.}, 34(3):553--579, 2005.

\bibitem[{Dujmovi\'{c} et~al.(2017)Dujmovi\'{c}, Morin, and
  Wood}]{dujmovic2017layered}
\textsc{Vida Dujmovi\'{c}, Pat Morin, and David~R. Wood}.
\newblock \href{https://doi.org/10.1016/j.jctb.2017.05.006}{Layered separators
  in minor-closed graph classes with applications}.
\newblock \emph{J. Combin. Theory Ser. B}, 127:111--147, 2017.

\bibitem[{Dujmovi\'{c} and Wood(2005)}]{DW2005subdivisions}
\textsc{Vida Dujmovi\'{c} and David~R. Wood}.
\newblock
  \href{http://www.dmtcs.org/volumes/abstracts/dm060221.abs.html}{Stacks,
  queues and tracks: layouts of graph subdivisions}.
\newblock \emph{Discrete Math. Theor. Comput. Sci.}, 7(1):155--201, 2005.

\bibitem[{Dvo{\v{r}}{\'a}k et~al.(2021)Dvo{\v{r}}{\'a}k, Huynh, Joret, Liu, and
  Wood}]{dvovrak2020notes}
\textsc{Zden{\v{e}}k Dvo{\v{r}}{\'a}k, Tony Huynh, Gwena{\"e}l Joret, Chun-Hung
  Liu, and David~R Wood}.
\newblock \href{https://doi.org/10.1007/978-3-030-62497-2_32}{Notes on graph
  product structure theory}.
\newblock In \textsc{David~R. Wood, Jan de~Gier, Cheryl~E. Praeger, and Terence
  Tao}, eds., \emph{2019-20 MATRIX Annals}, pp. 513--533. Springer, 2021.

\bibitem[{Dvo\v{r}\'{a}k(2013)}]{dvorak2014approximation}
\textsc{Zden\v{e}k Dvo\v{r}\'{a}k}.
\newblock \href{https://doi.org/10.1016/j.ejc.2012.12.004}{Constant-factor
  approximation of the domination number in sparse graphs}.
\newblock \emph{European J. Combin.}, 34(5):833--840, 2013.

\bibitem[{Dvo\v{r}\'{a}k and Norin(2016)}]{DN2016sublinear}
\textsc{Zden\v{e}k Dvo\v{r}\'{a}k and Sergey Norin}.
\newblock \href{https://doi.org/10.1137/15M1017569}{Strongly sublinear
  separators and polynomial expansion}.
\newblock \emph{SIAM J. Discrete Math.}, 30(2):1095--1101, 2016.

\bibitem[{Eppstein(2000)}]{eppstein2000diameter}
\textsc{David Eppstein}.
\newblock \href{https://doi.org/10.1007/s004530010020}{Diameter and treewidth
  in minor-closed graph families}.
\newblock \emph{Algorithmica}, 27(3-4):275--291, 2000.

\bibitem[{Eppstein(2015)}]{eppstein2015shallow}
\textsc{David Eppstein}.
\newblock
  \href{https://11011110.github.io/blog/2015/10/03/why-shallow-minors.html}{Why
  shallow minors matter for graph drawing --- 11011110}.
\newblock 2015.
\newblock [Online; accessed 15-October-2021].

\bibitem[{Eppstein and Gupta(2017)}]{EG2017road}
\textsc{David Eppstein and Siddharth Gupta}.
\newblock \href{https://doi.org/10.1145/3139958.3139999}{Crossing patterns in
  nonplanar road networks}.
\newblock In \emph{Proc. of the 25th ACM SIGSPATIAL International Conference on
  Advances in Geographic Information Systems}, SIGSPATIAL '17. ACM, 2017.

\bibitem[{Eppstein et~al.(2022)Eppstein, Hickingbotham, Merker, Norin, Seweryn,
  and Wood}]{EHMNSW22}
\textsc{David Eppstein, Robert Hickingbotham, Laura Merker, Sergey Norin,
  Micha\l{}~T. Seweryn, and David~R. Wood}.
\newblock \href{http://arxiv.org/abs/2202.05327}{Three-dimensional graph
  products with unbounded stack-number}.
\newblock 2022.
\newblock arXiv:2202.05327.

\bibitem[{Esperet et~al.(2020)Esperet, Joret, and Morin}]{EJM2020universal}
\textsc{Louis Esperet, Gwena{\"e}l Joret, and Pat Morin}.
\newblock \href{http://arxiv.org/abs/2010.05779}{Sparse universal graphs for
  planarity}.
\newblock 2020.
\newblock arXiv:2010.05779.

\bibitem[{Esperet and Raymond(2018)}]{ER2018expansion}
\textsc{Louis Esperet and Jean-Florent Raymond}.
\newblock \href{https://doi.org/10.1016/j.ejc.2017.09.003}{Polynomial expansion
  and sublinear separators}.
\newblock \emph{European J. Combin.}, 69:49--53, 2018.

\bibitem[{Esperet and Wiechert(2018)}]{EW2018boxicity}
\textsc{Louis Esperet and Veit Wiechert}.
\newblock \href{https://doi.org/10.37236/7787}{Boxicity, poset dimension, and
  excluded minors}.
\newblock \emph{Electron. J. Combin.}, 25(4):Paper No. 4.51, 11, 2018.

\bibitem[{Grohe et~al.(2018)Grohe, Kreutzer, Rabinovich, Siebertz, and
  Stavropoulos}]{GKRSS2018coveringsnowhere}
\textsc{Martin Grohe, Stephan Kreutzer, Roman Rabinovich, Sebastian Siebertz,
  and Konstantinos Stavropoulos}.
\newblock \href{https://doi.org/10.1137/18M1168753}{Coloring and covering
  nowhere dense graphs}.
\newblock \emph{SIAM J. Discrete Math.}, 32(4):2467--2481, 2018.

\bibitem[{Grohe et~al.(2017)Grohe, Kreutzer, and
  Siebertz}]{GKS2017propertiesnowhere}
\textsc{Martin Grohe, Stephan Kreutzer, and Sebastian Siebertz}.
\newblock \href{https://doi.org/10.1145/3051095}{Deciding first-order
  properties of nowhere dense graphs}.
\newblock \emph{J. ACM}, 64(3):Art. 17, 32, 2017.

\bibitem[{Hammack et~al.(2011)Hammack, Imrich, and
  Klav\v{z}ar}]{hammack2011handbook}
\textsc{Richard Hammack, Wilfried Imrich, and Sandi Klav\v{z}ar}.
\newblock Handbook of product graphs.
\newblock CRC Press, 2nd edn., 2011.

\bibitem[{Har-Peled and Quanrud(2017)}]{HQ2017lowdensity}
\textsc{Sariel Har-Peled and Kent Quanrud}.
\newblock \href{https://doi.org/10.1137/16M1079336}{Approximation algorithms
  for polynomial-expansion and low-density graphs}.
\newblock \emph{SIAM J. Comput.}, 46(6):1712--1744, 2017.

\bibitem[{Harvey and Wood(2017)}]{HW2017tied}
\textsc{Daniel~J. Harvey and David~R. Wood}.
\newblock \href{https://doi.org/10.1002/jgt.22030}{Parameters tied to
  treewidth}.
\newblock \emph{J. Graph Theory}, 84(4):364--385, 2017.

\bibitem[{Heath et~al.(1992)Heath, Leighton, and Rosenberg}]{HLR92}
\textsc{Lenwood~S. Heath, Frank~Thomson Leighton, and Arnold~L. Rosenberg}.
\newblock \href{https://doi.org/10.1137/0405031}{Comparing queues and stacks as
  machines for laying out graphs}.
\newblock \emph{SIAM J. Discrete Math.}, 5(3):398--412, 1992.

\bibitem[{Heath and Rosenberg(1992)}]{HR1992laying}
\textsc{Lenwood~S. Heath and Arnold~L. Rosenberg}.
\newblock \href{https://doi.org/10.1137/0221055}{Laying out graphs using
  queues}.
\newblock \emph{SIAM J. Comput.}, 21(5):927--958, 1992.

\bibitem[{Hickingbotham(2022)}]{H2022odd}
\textsc{Robert Hickingbotham}.
\newblock \href{http://arxiv.org/abs/2203.10402}{Odd colourings, conflict-free
  colourings and strong colouring numbers}.
\newblock 2022.
\newblock arXiv:2203.10402.

\bibitem[{Hickingbotham et~al.(2022)Hickingbotham, Jungeblut, Merker, and
  Wood}]{HJMW2022square}
\textsc{Robert Hickingbotham, Paul Jungeblut, Laura Merker, and David~R. Wood}.
\newblock \href{http://arxiv.org/abs/2203.03772}{The product structure of
  squaregraphs}.
\newblock 2022.
\newblock arXiv:2203.03772.

\bibitem[{Hickingbotham and Wood(2021)}]{HW2021products}
\textsc{Robert Hickingbotham and David~R. Wood}.
\newblock \href{http://arxiv.org/abs/2110.00721}{Structural properties of graph
  products}.
\newblock 2021.
\newblock arXiv:2110.00721.

\bibitem[{Hong and Tokuyama(2020)}]{HT2020beyond}
\textsc{Seok-Hee Hong and Takeshi Tokuyama}, eds.
\newblock \href{https://doi.org/10.1007/978-981-15-6533-5}{Beyond planar
  graphs}.
\newblock Springer, 2020, {C}ommunications of NII Shonan Meetings.

\bibitem[{Huynh et~al.(2021)Huynh, Mohar, {\v{S}}{\'a}mal, Thomassen, and
  Wood}]{HMSTW2021universality}
\textsc{Tony Huynh, Bojan Mohar, Robert {\v{S}}{\'a}mal, Carsten Thomassen, and
  David~R. Wood}.
\newblock \href{http://arxiv.org/abs/2109.00327}{Universality in minor-closed
  graph classes}.
\newblock 2021.
\newblock arXiv:2109.00327.

\bibitem[{Huynh and Wood(2021)}]{HW2021treedensities}
\textsc{Tony Huynh and David~R. Wood}.
\newblock \href{https://doi.org/10.4153/S0008414X21000316}{Tree densities in
  sparse graph classes}.
\newblock \emph{Canadian J. Math.}, 2021.

\bibitem[{Joret and Micek(2021)}]{JM2021weak}
\textsc{Gwena{\"e}l Joret and Piotr Micek}.
\newblock \href{http://arxiv.org/abs/2102.10061}{Improved bounds for weak
  coloring numbers}.
\newblock 2021.
\newblock arXiv:2102.10061.

\bibitem[{Kaufmann and Ueckerdt(2022)}]{KU2014density}
\textsc{Michael Kaufmann and Torsten Ueckerdt}.
\newblock \href{https://doi.org/https://doi.org/10.37236/10521}{The density of
  fan-planar graphs}.
\newblock \emph{Electron. J. Combin.}, 29(1):Paper No. 1.29, 2022.

\bibitem[{Kierstead and Trotter(1994)}]{KT1994uncooperative}
\textsc{Hal~A. Kierstead and William~T. Trotter}.
\newblock \href{https://doi.org/10.1002/jgt.3190180605}{Planar graph coloring
  with an uncooperative partner}.
\newblock \emph{J. Graph Theory}, 18(6):569--584, 1994.

\bibitem[{Kierstead and Yang(2003)}]{KY2003orderings}
\textsc{Hal~A. Kierstead and Daqing Yang}.
\newblock \href{https://doi.org/10.1023/B:ORDE.0000026489.93166.cb}{Orderings
  on graphs and game coloring number}.
\newblock \emph{Order}, 20(3):255--264, 2003.

\bibitem[{Klemz et~al.(2021)Klemz, Knorr, Reddy, and
  Schr{\"{o}}der}]{KKRS2021fan}
\textsc{Boris Klemz, Kristin Knorr, Meghana~M. Reddy, and Felix
  Schr{\"{o}}der}.
\newblock \href{https://doi.org/10.1007/978-3-030-92931-2\_4}{Simplifying
  non-simple fan-planar drawings}.
\newblock In \textsc{Helen~C. Purchase and Ignaz Rutter}, eds., \emph{Proc.
  29th International Symposium on Graph Drawing and Network Visualization, {GD}
  2021}, vol. 12868 of \emph{Lecture Notes in Computer Science}, pp. 57--71.
  Springer, 2021.

\bibitem[{Kostochka et~al.(1997)Kostochka, Sopena, and Zhu}]{KSZ1997acyclic}
\textsc{Alexandr~V. Kostochka, Eric Sopena, and Xuding Zhu}.
\newblock
  \href{https://doi.org/10.1002/(SICI)1097-0118(199704)24:4<331::AID-JGT5>3.0.CO;2-P}{Acyclic
  and oriented chromatic numbers of graphs}.
\newblock \emph{J. Graph Theory}, 24(4):331--340, 1997.

\bibitem[{K\"{u}ndgen and Pelsmajer(2008)}]{KP2008nonrepetitive}
\textsc{Andr\'{e} K\"{u}ndgen and Michael~J. Pelsmajer}.
\newblock \href{https://doi.org/10.1016/j.disc.2007.08.043}{Nonrepetitive
  colorings of graphs of bounded tree-width}.
\newblock \emph{Discrete Math.}, 308(19):4473--4478, 2008.

\bibitem[{Leighton(1983)}]{leighton1983VLSI}
\textsc{Frank~Thomson Leighton}.
\newblock Complexity issues in {VLSI}: Optimal layouts for the shuffle-exchange
  graph and other networks.
\newblock MIT Press, 1983.

\bibitem[{Mohar and Thomassen(2001)}]{MT2001surfaces}
\textsc{Bojan Mohar and Carsten Thomassen}.
\newblock Graphs on surfaces.
\newblock Johns Hopkins University Press, 2001.

\bibitem[{Ne{\v{s}}et{\v{r}}il and Ossona~de
  Mendez(2012)}]{nevsetvril2012sparsity}
\textsc{Jaroslav Ne{\v{s}}et{\v{r}}il and Patrice Ossona~de Mendez}.
\newblock Sparsity: graphs, structures, and algorithms, vol.~28.
\newblock Springer, 2012.

\bibitem[{Ne\v{s}et\v{r}il and Ossona~de Mendez(2008)}]{nevsetvril2008grad}
\textsc{Jaroslav Ne\v{s}et\v{r}il and Patrice Ossona~de Mendez}.
\newblock \href{https://doi.org/10.1016/j.ejc.2006.07.013}{Grad and classes
  with bounded expansion. {I}. {D}ecompositions}.
\newblock \emph{European J. Combin.}, 29(3):760--776, 2008.

\bibitem[{Ne\v{s}et\v{r}il et~al.(2012)Ne\v{s}et\v{r}il, Ossona~de Mendez, and
  Wood}]{NOW2012examples}
\textsc{Jaroslav Ne\v{s}et\v{r}il, Patrice Ossona~de Mendez, and David~R.
  Wood}.
\newblock \href{https://doi.org/10.1016/j.ejc.2011.09.008}{Characterisations
  and examples of graph classes with bounded expansion}.
\newblock \emph{European J. Combin.}, 33(3):350--373, 2012.

\bibitem[{Ollmann(1973)}]{ollmann1973book}
\textsc{L.~Taylor Ollmann}.
\newblock On the book thicknesses of various graphs.
\newblock In \textsc{Frederick Hoffman, Roy~B. Levow, and Robert S.~D. Thomas},
  eds., \emph{Proc. 4th Southeastern Conference on Combinatorics, Graph Theory
  and Computing}, vol. VIII of \emph{Congr. Numer.}, p. 459. 1973.

\bibitem[{Ossona~de Mendez et~al.(2019)Ossona~de Mendez, Oum, and
  Wood}]{OOW2019defective}
\textsc{Patrice Ossona~de Mendez, Sang-Il Oum, and David~R. Wood}.
\newblock \href{https://doi.org/10.1007/s00493-018-3733-1}{Defective colouring
  of graphs excluding a subgraph or minor}.
\newblock \emph{Combinatorica}, 39(2):377--410, 2019.

\bibitem[{Pach and T\'{o}th(1997)}]{PT1997crossings}
\textsc{J\'{a}nos Pach and G\'{e}za T\'{o}th}.
\newblock \href{https://doi.org/10.1007/BF01215922}{Graphs drawn with few
  crossings per edge}.
\newblock \emph{Combinatorica}, 17(3):427--439, 1997.

\bibitem[{Pilipczuk and Siebertz(2021)}]{PS2019centred}
\textsc{Micha\l Pilipczuk and Sebastian Siebertz}.
\newblock \href{https://doi.org/10.1016/j.jctb.2021.06.002}{Polynomial bounds
  for centered colorings on proper minor-closed graph classes}.
\newblock \emph{J. Combin. Theory Ser. B}, 151:111--147, 2021.

\bibitem[{Pupyrev(2020)}]{pupyrev2020book}
\textsc{Sergey Pupyrev}.
\newblock \href{http://arxiv.org/abs/2007.15102}{Book embeddings of graph
  products}.
\newblock 2020.
\newblock arXiv:2007.15102.

\bibitem[{Reed(1997)}]{reed1997treewidth}
\textsc{Bruce Reed}.
\newblock \href{https://doi.org/10.1017/CBO9780511662119.006}{Tree width and
  tangles: a new connectivity measure and some applications}.
\newblock In \emph{Surveys in combinatorics, 1997}, vol. 241 of \emph{London
  Math. Soc. Lecture Note Ser.}, pp. 87--162. Cambridge Univ. Press, 1997.

\bibitem[{Robertson and Seymour(1986)}]{robertson1986algorithmic}
\textsc{Neil Robertson and Paul Seymour}.
\newblock \href{https://doi.org/10.1016/0196-6774(86)90023-4}{Graph minors.
  {II}. {A}lgorithmic aspects of tree-width}.
\newblock \emph{J. Algorithms}, 7(3):309--322, 1986.

\bibitem[{Robertson and Seymour(2003)}]{robertson2003graph}
\textsc{Neil Robertson and Paul Seymour}.
\newblock \href{https://doi.org/10.1016/S0095-8956(03)00042-X}{Graph minors.
  {XVI}. {E}xcluding a non-planar graph}.
\newblock \emph{J. Combin. Theory Ser. B}, 89(1):43--76, 2003.

\bibitem[{Scott and Wood(2020)}]{SW2020boxicity}
\textsc{Alex Scott and David~R. Wood}.
\newblock \href{https://doi.org/10.1090/tran/7962}{Better bounds for poset
  dimension and boxicity}.
\newblock \emph{Trans. Amer. Math. Soc.}, 373(3):2157--2172, 2020.

\bibitem[{Thomassen(1986)}]{thomassen1986interval}
\textsc{Carsten Thomassen}.
\newblock \href{https://doi.org/10.1016/0095-8956(86)90061-4}{Interval
  representations of planar graphs}.
\newblock \emph{J. Combin. Theory Ser. B}, 40(1):9--20, 1986.

\bibitem[{Ueckerdt et~al.(2021)Ueckerdt, Wood, and Yi}]{UWY2021GPST}
\textsc{Torsten Ueckerdt, David~R Wood, and Wendy Yi}.
\newblock \href{http://arxiv.org/abs/2108.00198}{An improved planar graph
  product structure theorem}.
\newblock 2021.
\newblock arXiv:2108.00198.

\bibitem[{van~den Heuvel et~al.(2017)van~den Heuvel, Ossona~de Mendez, Quiroz,
  Rabinovich, and Siebertz}]{HMQRS2017fixed}
\textsc{Jan van~den Heuvel, Patrice Ossona~de Mendez, Daniel Quiroz, Roman
  Rabinovich, and Sebastian Siebertz}.
\newblock \href{https://doi.org/10.1016/j.ejc.2017.06.019}{On the generalised
  colouring numbers of graphs that exclude a fixed minor}.
\newblock \emph{European J. Combin.}, 66:129--144, 2017.

\bibitem[{van~den Heuvel and Wood(2018)}]{HW2018improper}
\textsc{Jan van~den Heuvel and David~R. Wood}.
\newblock \href{https://doi.org/10.1112/jlms.12127}{Improper colourings
  inspired by {H}adwiger's conjecture}.
\newblock \emph{J. Lond. Math. Soc. (2)}, 98(1):129--148, 2018.
\newblock arXiv:1704.06536.

\bibitem[{Wiechert(2017)}]{wiechert2017queue}
\textsc{Veit Wiechert}.
\newblock \href{https://doi.org/https://doi.org/10.37236/6429}{On the
  queue-number of graphs with bounded tree-width}.
\newblock \emph{Electron. J. Combin.}, 24(1):Paper No. 1.65, 18, 2017.

\bibitem[{Wood(2005)}]{wood2005queue}
\textsc{David~R. Wood}.
\newblock \href{http://dmtcs.episciences.org/352}{Queue layouts of graph
  products and powers}.
\newblock \emph{Discrete Math. Theor. Comput. Sci.}, 7(1):255--268, 2005.

\bibitem[{Wood(2021)}]{wood2020nonrepetitive}
\textsc{David~R. Wood}.
\newblock \href{https://doi.org/10.37236/9777}{Nonrepetitive graph colouring}.
\newblock \emph{Electron. J. Combin}, (DS24), 2021.

\bibitem[{Zhu(2009)}]{zhu2009generalized}
\textsc{Xuding Zhu}.
\newblock \href{https://doi.org/10.1016/j.disc.2008.03.024}{Colouring graphs
  with bounded generalized colouring number}.
\newblock \emph{Discrete Math.}, 309(18):5562--5568, 2009.

\end{thebibliography}
\end{document}